\documentclass[12pt]{article}

\usepackage[top=1in,bottom=1in, left=1in, right=1in]{geometry}

\usepackage{amsmath,amsthm,amssymb,amsopn,amsfonts,pdfpages,algorithmic,algorithm,stmaryrd}
\usepackage{graphics,dsfont}
\usepackage{graphicx}
\usepackage{bbm,subfig}
\usepackage{caption,url}
\usepackage{authblk}
\usepackage[colorlinks=true,linkcolor=cyan,citecolor=blue]{hyperref}

\usepackage{todonotes}

\newcommand{\pnc}{\mathcal{P}_n^\mathcal{C}}

\newcommand{\twospace}{\renewcommand{\baselinestretch}{1.3}\normalsize}

\renewcommand{\Re}{\mathbb{R}}
\newcommand{\argmin}{\text{argmin}}
\newcommand{\argmax}{\text{argmax}}

\renewcommand{\Pr}{\mathbb{P}}
\newcommand{\tA}{\widetilde A}
\newcommand{\tB}{\widetilde B}

\newcommand{\ck}{\mathcal{K}}
\newcommand{\hA}{\widehat A}
\newcommand{\hB}{\widehat B}
\newcommand{\hQ}{\widehat Q}

\newcommand{\p}{\mathbb{P}}
\newcommand{\e}{\mathbb{E}}

\newcommand{\td}{\text{oracle}}
\newcommand{\hd}{\widehat\delta}
\newcommand{\gd}{\delta^{\text{GMP}}}
\newcommand{\tb}{\textcolor{black}}

\newcommand{\Ex}{\mathbb{E}}

\newcommand{\HER}{\text{CorrER}(Q_1,Q_2,R)}

\newcommand{\id}{\mathrm{id}}

\newtheorem{theorem}{Theorem}
\newtheorem{proposition}[theorem]{Proposition}
\newtheorem{lemma}[theorem]{Lemma}

\theoremstyle{definition}
\newtheorem{definition}[theorem]{Definition}
\theoremstyle{remark}
\newtheorem{remark}[theorem]{Remark}
\newtheorem{example}[theorem]{Example}

\twospace
\begin{document}
\title{Matchability of heterogeneous networks pairs}
\author[$\dag$]{Vince~Lyzinski}
\author[$\ddag$]{Daniel L. Sussman}

\affil[$\dag$]{\small Department of Mathematics and Statistics, University of Massachusetts Amherst
   }
\affil[$\ddag$]{\small Department of Mathematics and Statistics, Boston University}

\maketitle

\begin{abstract}
We consider the problem of graph matchability in non-identically distributed networks.
In a general class of edge-independent networks, we demonstrate that graph matchability can be lost with high probability when matching the networks directly.
We further demonstrate that under mild model assumptions, matchability is almost perfectly recovered by centering the networks using Universal Singular Value Thresholding before matching.
These theoretical results are then demonstrated in both real and synthetic simulation settings.
We also recover analogous core-matchability results in a very general core-junk network model, wherein some vertices do not correspond between the graph pair.
\end{abstract}

\section{Introduction and Background}

The graph matching problem seeks to find an alignment between the vertex sets of two graphs that best preserves common structure across graphs. 
At its simplest, it can be formulated as follows:  Given two adjacency matrices $A$ and $B$ corresponding to $n$-vertex graphs, the \textit{graph matching problem} seeks to minimize 
$\|A-PBP^T\|_F$ over permutation matrices $P$; i.e., the graph matching problem seeks a relabeling of the vertices of $B$ that minimizes the number of induced edge disagreements between $A$ and $PBP^T$.
Variants and extensions of this problem have been extensively studied in the literature, with applications across areas as diverse as biology and neuroscience \cite{zaslavskiy2009,klau,chen2016joint,prot}, computer vision \cite{robles, escolano,lin2010layered}, pattern recognition \cite{lee2010graph,sang2012robust,zhou}, and social network analysis \cite{limatching,zhang2016final,regal}, among others.
For a survey of many of the recent applications and approaches to the graph matching problem, see the sequence of survey papers \cite{ConteReview,foggia2014graph,Emmert-Streib2016-st}. 
While recent results \cite{babai2016graph} have whittled away at the complexity of the related \emph{graph isomorphism problem}---determining whether a permutation matrix $P$ exists satisfying $A=PBP^T$---
at its most general, where $A$ and $B$ are allowed to be weighted and directed, the graph matching problem is known to be NP-hard.
Indeed, in this case, the graph matching problem is equivalent to the notoriously difficult quadratic assignment problem \cite{loiola2007survey,cela2013quadratic,qapref}.
However, recent approaches that leverage efficient representation/learning methodologies (see, for example, \cite{netalign,zhang2016final,regal}) have shown excellent empirical performance matching networks with up to millions of nodes.

In addition to algorithmic advancements in graph matching, there has been a flurry of activity studying the closely related problem of graph matchability:  Given a latent alignment between the vertex sets of two graphs, can graph matching uncover this alignment in the presence of shuffled vertex labels?
This problem arises in a variety of contexts, from network de-anonymization and privatization to multi-network hypothesis testing \cite{lyzinski2016information} to multimodality graph embedding methodologies \cite{jointLi}.
Many existing results are concerned with recovering a latent alignment present across random graph models where each of $A$ and $B$ have identical marginal distributions, and exciting advancements on the threshold of matchable versus unmatchable graphs have been made across many random graph settings, including: the homogeneous correlated Erd\H os-Renyi model (see, for example, \cite{pedarsani2011privacy,JMLR:v15:lyzinski14a,eralg,eralg2}), the correlated stochastic blockmodel setting (see, for example, \cite{onaran2016optimal,lyzinski2016information}), the $\rho$-correlated heterogeneous Erd\H os-Renyi model (see, for example, \cite{rel,vncon}), and in the correlated heterogeneous Erd\H os-Renyi model with varying edge correlations (see, for example, \cite{sussman2018matched,lyzinski2017consistent}).
In the non-identically distributed model setting, the work in \cite{het2,het1,het3} provide theoretic phase transitions on matchability in the $(A,B)\sim$Erd\H os-R\'enyi(p,q,$\varrho$) model (i.e., $A\sim$Erd\H os-R\'enyi(n,p), $B\sim$Erd\H os-R\'enyi(n,q) and the edge correlation across graphs in provided by the constant $\varrho$; see Definition \ref{def:corrER}).

\tb{The above results range from providing theoretic phase transitions on matchability \cite{het2,het1,lyzinski2016information} to providing nearly efficient methods for achieving matchability from an algorithmic perspective \cite{eralg2,eralg,het3,fang2018tractable}.
While they have served to establish a novel theoretical understanding of the matchability problem, in each case the transition from matchable to unmatchable graphs is defined in terms of decreasing across-graph correlation and within-graph sparsity.
Importantly, it is not a function of fundamentally different probabilistic structures across the graphs to be matched.
As we often witness in applications, the graph topologies can differ significantly even amongst vertices that correspond to the same entity across networks.}
Social networks offer a compelling example of this, where matching across different social network platforms requires the understanding that not all users will be behaving homogeneously across different network platforms \cite{limatching}.
Both theoretically (see Example \ref{ex:bad}) and practically (see, for example, \cite{chen2016joint}), this distributional heterogeneity can have a deleterious effect on graph matchability.

Herein, we propose one possible solution for ameliorating the effect of $A$ and $B$ not being identically distributed, namely via a Universal Singular Value Thresholding \cite{chatterjee2015matrix} (USVT) centering preprocessing step.
Working in a general correlated edge-independent random graph model (see Definition \ref{def:GMP}), we theoretically demonstrate that USVT centering asymptotically almost surely recovers the matchability for all but a vanishing fraction of the nodes (see Theorem \ref{thm:estimatedcenteredmatched}).
In addition, we recover analogous result (see Theorem \ref{thm:coreestcenteredmatched}) in the setting in which only a fraction of the vertices possess a true correspondence across networks, generalizing and extending the results of \cite{pedarsani2011privacy,kazemi2015can,yartseva2013performance}.

This centering step is practically implementable on even very large networks and is demonstrated to have a significant positive impact on graph matchability in both real and synthetic data settings (see Section \ref{sec:expsim}).
While the results contained herein do not guarantee that any computationally efficient algorithm will be able to perfectly (or almost perfectly) align any given networks after USVT centering, they provide a theoretical foundation for subsequently studying algorithmic effectiveness.
Indeed, they ensure that with high probability the optimal alignment according to the graph matching objective function is, essentially, the true latent vertex alignment, guaranteeing that subsequent optimization procedures are, at the least, seeking the right permutation.

\subsection{Notation}
\tb{The following notation will be used throughout the manuscript:  for $k\in\mathbb{Z}>0$, we will let $J_k$ denote the hollow $k\times k$ matrix with all $1$'s on its off-diagonal, $[0]_k$ will denote the $k\times k$ matrix of all $0$'s, and $[k]$ will denote the set $\{1,2,\ldots,k\}$.
We will consider $A$ and $B$ interchangeably as adjacency matrices and as the corresponding graphs consisting of vertices and edges.
For a set $S\subset V(A)$, we denote by $A[S]$ the induced subgraph of $A$ on the vertices of $S$.}

\tb{For a matrix $A\in\mathbb{R}^{n\times n}$, the Froebenius norm of the matrix is defined as 
$$\|A\|_F=\sqrt{\sum_{i=1}^n \sum_{j=1}^n A(i,j)^2 },$$
and 
the operator norm of $A$ is denoted 
$$\|A\|=\sigma_{1 }(A),$$
where $\sigma_{1 }(A)$ is the largest singular value of $A$.
We denote the $\text{trace}$ of $A$ via
$$\text{tr}(A)=\sum_{i=1}^n A(i,i).$$
Below we will make use of the following trace form of the Frobenius norm: $\|A\|_F^2=\text{tr}(A^TA)$; see \cite{horn2012matrix} for more on the Frobenius norm and its many uses.
For matrices $A\in\mathbb{R}^{n\times n}$ and $B\in\mathbb{R}^{m\times m}$, we define $A\oplus B\in\in\mathbb{R}^{(n+m)\times (n+m)}$ via 
$$
A\oplus B=\begin{pmatrix}
A & [0]_{n,m}\\
[0]_{m,n} & B
\end{pmatrix},
$$
where $[0]_{n,m}$ is the $n\times m$ matrix of all $0$'s.}

We will also make extensive use of modern asymptotic notation.
To review, if $f$ and $g$ are nonnegative functions of $n\in \mathbb{Z}\geq 0$, then we write 
\begin{align*}
f(n)=o(g(n))	&\text{ if }\lim _{n\rightarrow \infty} \frac {f(n)}{g(n)}=0;\\
f(n)=O(g(n)) &\text{ if }\exists\, C>0, \exists\, n_0>0\text{ s.t. }f(n)<C g(n)\text{ for all }n\geq n_0;\\
f(n)=\omega(g(n)) &\text{ if } g(n)=o(f(n));\\
f(n)=\Omega(g(n)) &\text{ if }g(n)=O(f(n));\\
f(n)=\Theta (g(n))&\text{ if }	f(n)=O(g(n))\text{ and }f(n)=\Omega(g(n)).
\end{align*}


\subsection[Correlated heterogeneous Erdos-Renyi graphs]{Correlated heterogeneous Erd\H os-R\'enyi graphs}
\label{sec:RcorrER}

Formally the graph matching problem (GMP) we will consider is defined as follows.
\begin{definition}
\label{def:GMP}
Let $A,B\in\mathbb{R}^{n\times n}$ be the adjacency matrices of weighted, undirected graphs on $n$ vertices. 
The graph matching problem is to find an element of
\begin{equation}
\label{eq:GMP}
\argmin_{P\in\Pi(n)}\|A-PBP^T\|_F^2=\argmin_{P\in\Pi(n)}-\mathrm{tr}(APBP^T),
\end{equation}
where $\Pi(n)$ is the set of $n \times n$ permutation matrices.
\end{definition}
\noindent Eq. (\ref{eq:GMP}) follows here from 
\begin{align*}
\|A-PBP^T\|_F^2&=\text{tr}\left((A-PBP^T)^T(A-PBP^T)\right)\\
&=\text{tr}(A^TA)-2\text{tr}(A^TPBP^T)+\text{tr}\left((PBP^T)^T(PBP^T)\right)\\
&=\text{tr}(A^TA)-2\text{tr}(APBP^T)+\text{tr}\left(B^TB\right)
\end{align*}

\tb{We note here that, traditionally, the GMP formulated in Definition \ref{def:GMP} is defined for \emph{unweighted} graphs $A$ and $B$.  
The extension we consider to weighted graphs is commonly used in the literature (see, for example, the work in \cite{umeyama1988eigendecomposition}) and is useful for studying situations in which edges/vertices in the network have weight features attached to them.
This added flexibility will be needed for subsequent theoretical developments and data applications.}

In the presence of a latent vertex alignment, $\phi:[n]\mapsto [n]$, between the vertices of $A$ and $B$, we wish to understand the extent to which graph matching $A$ and $B$ will recover $\phi$; i.e., if $P_\phi$ is the permutation matrix corresponding to $\phi$, will $\{P_\phi\}=\text{argmin}_{P\in\Pi(n)}\|A-PBP^T\|_F^2$?
In order to study this problem from a probabilistic perspective, we introduce a bivariate random graph model with a natural vertex alignment across graphs: the bivariate, correlated, heterogeneous, Erd\H os-R\'enyi random graph.  
\begin{definition}
\label{def:corrER}
For $R,Q_1,Q_2\in\mathbb{R}^{n\times n}$ symmetric, matrices, we say $A,B$ are instantiations of the $R$-correlated Heterogeneous Erd\H os-R\'enyi random graph model with parameters $R,Q_1,Q_2$  (abbreviated as $\HER$) if:
\begin{itemize}
\item[1.]  $A\sim$ER$(Q_1)$; i.e., $A$ is an independent edge random graph with no self-loops satisfying
$$\mathds{1}_{\{u,v\}\in E(A)}\sim \text{Bernoulli}(Q_1(u,v))$$ for each $\{u,v\}\in\binom{V}{2};$ 
\item[2.]  $B\sim$ER$(Q_2)$; i.e., $B$ is an independent edge random graph with no self-loops satisfying
$$\mathds{1}_{\{u,v\}\in E(B)}\sim \text{Bernoulli}(Q_2(u,v))$$ for each $\{u,v\}\in\binom{V}{2}$;
\item[3.] Edges across networks are collectively independent except that for each $\{u,v\}\in \binom{V}{2}$, the correlation between $A(u,v)$ and $B(u,v)$ is 
$$R(u,v)=\frac{\text{Cov}(A(u,v),B(u,v))}{\sqrt{\text{Var}(A(u,v))\text{Var}(B(u,v)) }}\geq0$$
\end{itemize}
 \end{definition}
Before proceeding further, we will make a few remarks on the $\HER$ random graph model.  
In the homogeneous ER$(n,p)$ model, network growth as $n\rightarrow\infty$ is natural, and we can consider an asymptotic regime in which $p$ depends on $n$.
Here, we similarly consider $Q_1$ and $Q_2$ to be dependent on $n$, but make no further assumptions on expressly how the dependence on $n$ is manifest.
This allows for us to consider classical homogeneous Erd\H os-R\'enyi, stochastic blockmodels \cite{sbm}, random dot product graphs (conditioned on the latent positions) \cite{young2007random}, etc., as subfamilies of our random graph model.

In addition, by allowing $Q_1$ and $Q_2$ to differ, this model allows for a latent correspondence to exist in settings where the underlying topology and degree structure of the graphs to be matched differs significantly.    
This distributional heterogeneity is often observed in real data settings (see, for example, the connectomes being aligned in \cite{chen2016joint} and the social networks aligned in \cite{limatching}), and we seek to understand the limitations of graph matching approaches when attempting to overcome this heterogeneity.
Note also that when $Q_1\neq Q_2$, there are restrictions on feasible correlations $R$: Indeed, if $X\sim$Bernoulli$(p)$ and $Y\sim$Bernoulli$(q)$ are $\varrho$-correlated with $p\geq q$, then the correlation must satisfy $\varrho\leq \sqrt{\frac{q(1-p)}{p(1-q)}}$.

Lastly, this model naturally allows us to consider a partition of $V$ into core ($\mathcal{C}$) and junk ($\mathcal{J}$) vertices $V=\mathcal{C}\cup\mathcal{J}$;
core vertices are those that have a corresponding vertex (i.e., true match) across networks, while junk vertices do not.
If we consider $(A,B)\sim\HER$ with $R$ of the form $R=R_c\oplus 0_{n_j}$ where $R_c\in(0,1]^{n_c\times n_c}$, then it is reasonable to define $\mathcal{C}=[n_c]$ and $\mathcal{J}=[n]\setminus[n_c]$.
For all $u\in\mathcal{J}$ it would then hold that $\text{Corr}(A(u,v),B(u,v))=0$ for all $v$, and $A(u,v)$ and $B(u,v)$ are independent random variables.   
A natural question to ask is when an optimal GM algorithm will correctly align the vertices in $\mathcal{C}$ across networks.
This problem was studied in the context of homogeneous ER networks with constant correlation in \cite{kazemi2015can}, and the results in Section \ref{sec:cores} generalize and extend those in \cite{kazemi2015can} to this more adaptable network model.

\begin{remark}
In what follows, $Q_1$ and $Q_2$ are not necessarily assumed to be hollow matrices.
We do assume our graphs are loop-free, so that $\e(A)$ (resp., $\e(B)$) are necessarily hollow and need not equal $Q_1$ (resp., $Q_2$). 
\end{remark}
\section{Graph matchability}

In the $A,B\sim\HER$ setting, we seek to understand when a graph matching procedure could correctly align the vertices across networks; \tb{i.e., if $\argmin_{P\in \Pi(n)} \|A-PBP^T\|_F=\{I_n\},$
 where $I_n$ denotes the identity matrix.
More generally, if $\delta(\cdot,\cdot):\mathbb{R}^{n\times n}\times \mathbb{R}^{n\times n}\mapsto \mathbb{R}$ is a dissimilarity (i.e., if it is a symmetric, non-negative function with $\delta(X,X)=0$ for all $X\in\mathbb{R}^{n\times n}$; see, for example, \cite{trosset2008semisupervised}), when is it the case that 
$$\argmin_{P\in \Pi(n)} \delta(A,PBP^T)=\{I_n\}?$$
In this more general framework, we consider the following definition of graph matchability:}
\begin{definition}
\label{def:matchable}
Let $\delta(\cdot,\cdot):\mathbb{R}^{n\times n}\times \mathbb{R}^{n\times n}\mapsto \mathbb{R}$ be a dissimilarity.
We will say that $(A,B)\sim \HER$ are $\delta$-matchable if 
$$\argmin_{P\in \Pi(n)} \delta(A,PBP^T)=\{I_n\},$$
 where $I_n$ denotes the identity matrix.
\end{definition}

\tb{By considering an appropriate $\delta$ in Definition \ref{def:matchable}, we can fit the classical GMP in the formulation; indeed, the GMP of Definition \ref{def:GMP} considers the dissimilarity defined via 
$$\delta(A,PBP^T)=\|A-PBP^T\|_F.$$
In this paper, we will consider, more generally, dissimilarity functions of the form 
$$\delta(A,PBP^T)=\|f(A)-f(PBP^T)\|_F,$$
for suitably defined matrix-valued function
\begin{align*}
f(\cdot)&:\mathbb{R}^{n\times n}\mapsto \mathbb{R}^{n\times n};
\end{align*}
Special cases of interest in our present $(A,B)\sim \HER$ setting are 
\begin{align}
\label{eq:delid}
\delta^{\text{GMP}}(A,PBP^T)&=\|f(A)-f(PBP^T)\|_F,\text{ where }f(x)=x;\\
\label{eq:delhat}
\widehat\delta(A,PBP^T)&=\|f(A)-f(PBP^T)\|_F,\text{ where }f(x)=x-\hQ_{x},
\end{align}
where $\hQ_{A}$ (resp., $\hQ_{B}$) is a suitable estimate of $Q_1$ (resp., $Q_2$) derived from $A$ (resp., $B$).
In addition to the notion of $\delta$-matchability for a dissimilarity $\delta$, we will also define the notion of \emph{oracle}-matchability.
We will say that $A$ and $B$ distributed $(A,B)\sim\HER$ are $\td$-matchable if
$$\argmin_{P\in\Pi(n)}\|A-\e(A)-P(B-\e(B))P^T\|_F=\{I_n\}.$$
In the sequel, oracle-matchability will provide a useful theoretic bridge between $\gd$-matchability and $\hd$-matchability. 
Note that we will write $\gd$-matchability and $\hd$-matchability for the notions defined in Eq. (\ref{eq:delid}) and (\ref{eq:delhat}) respectively.}

\tb{A natural question to ask is why we define the GMP in terms of $\gd$ and not in terms of a more general dissimilarity $\delta$; indeed, alternate dissimilarities $\delta$ have been considered in the definition of the GMP in the graph matching literature (see for example \cite{zhang2018consistent,zhang2018unseeded}).
Moreover, we consider the GMP objective function formulation in Definition \ref{def:GMP} even though in numerous settings the optimal solution to this GMP may not be a given latent vertex alignment, and in this section, we will see instances of when $A,B\sim\HER$ are, with high probability, not $\gd$-matchable.
Our choice of $\gd$ in the GMP is motivated by two main factors.
First, this is the classical definition of graph matching and ties our current work to a vast graph matching literature.
Second, we seek to understand conditions for when the original $\gd$-matchability fails, yet there is a suitable dissimilarity $\delta$ for which $\delta$-matchability is achieved.
As the formulation in Definition \ref{def:GMP} is commonly used in practice, this could provide practical guidance for when vertex labels can be recovered via a different objective function viewpoint.}

\tb{In recent work addressing the question of $\gd$-matchability, results have been established for the $R=\varrho J_n$, $Q_1=Q_2=pJ_n$ setting 
(see, for example, \cite{pedarsani2011privacy,JMLR:v15:lyzinski14a,eralg,eralg2}), in the correlated stochastic blockmodel setting (see, for example, \cite{onaran2016optimal,lyzinski2016information}), in the correlated heterogeneous Erd\H os-Renyi model (see, for example, \cite{rel,vncon}), and in the general $R$ and general $Q_1=Q_2=Q$ setting (see, for example, \cite{sussman2018matched,lyzinski2017consistent}).
In the non-identically distributed model setting, the work in \cite{het2,het1,het3} considers $R=\varrho J_n$, $Q_1=pJ_n$, and $Q_2=qJ_n$.
In each setting, the results showed that for sufficiently dense, sufficiently correlated graphs, $\gd$-matchability is almost surely achieved.
Converse results in \cite{het2,het1,lyzinski2016information,JMLR:v15:lyzinski14a} show that in the sufficiently sparse and/or weakly correlated setting, $\gd$-matchability is a.s. lost (i.e., a.s. the solution to the GMP is not the latent alignment).
The work in \cite{het2,het1} deserves special mention, as the converse results therein are proven  for general $\delta$-matchability; i.e., for sufficiently sparse and/or weakly correlated networks, $\delta$-matchability is a.s.\@ lost for all dissimilarities $\delta$. }

\tb{In these examples, it is sparsity and/or weak dependence that is potentially thwarting the matching in each instance and not the heterogeneity of the model itself.
As the next straightforward} but illustrative example demonstrates, the degree and structural heterogeneity across networks allowed for in the $\HER$ model makes the question of $\gd$-matchability a bit more nuanced.

\begin{example}
\label{ex:bad}
Consider the following correlated heterogeneous stochastic blockmodel example.
Let $p,q,r\in[0,1]$ be distinct, and define
$$Q_1=\begin{pmatrix}
pJ_n& rJ_n\\
rJ_n& qJ_n
\end{pmatrix}
\text{ and }
Q_2=\begin{pmatrix}
qJ_n& rJ_n\\
rJ_n& pJ_n
\end{pmatrix}.$$
Let $V_1=[n]$ be the vertices in block 1 in $A$ and $B$ and let $V_2=[2n]\setminus[n]$ be the vertices in block 2.  
Assuming $p>q$, and letting $\varrho_1\leq\sqrt{\frac{q(1-p)}{p(1-q)}}$ and $\varrho_2\in[0,1]$ we consider $R$ of the form
$$R=\begin{pmatrix}
\varrho_1J_n& \varrho_2J_n\\
\varrho_2J_n& \varrho_1J_n
\end{pmatrix}.$$
\tb{Unlike in the cases where the loss of $\gd$-matchability is due to network sparsity and/or weak correlation, in this example the non-identically distributed nature of $A$ and $B$ can obfuscate the true alignment from a graph matching perspective.
Indeed, for many choices of the parameters above, the optimal permutation for the GMP in Definition \ref{def:GMP} will \text{not} be the latent correspondence, and permuting blocks 1 and 2 will, with high probability, yield a better GMP objective function.}

\tb{To wit, let $P\in\Pi(2n)$ be any permutation such that $PQ_2P^T=Q_1$ (so that $P$ aligns block 1 in $A$ to block $2$ in $B$ and vice versa) with corresponding permutation $\tau$.
The number of edges $\{u,v\}$ such that $\{\tau(u),\tau(v\}=\{u,v\}$ is bounded above by $n$ and 
\begin{align*}
\left\{  \{u,v\}\in\binom{V_1}{2}\right\}&\cap \bigg\{  \{u,v\}\text{ s.t. }\{\tau(u),\tau(v\}=\{u,v\}\bigg\}=\emptyset\\
\left\{  \{u,v\}\in\binom{V_2}{2}\right\}&\cap \bigg\{  \{u,v\}\text{ s.t. }\{\tau(u),\tau(v\}=\{u,v\}\bigg\}=\emptyset.
\end{align*}
Therefore, 
\begin{align*}
\frac{1}{2}\e\text{tr}(APBP^T)&\geq \binom{n}{2}(p^2+q^2)+(n^2-n)r^2\\
\frac{1}{2}\e\text{tr}(AB)&= 2\binom{n}{2}(\varrho_1\sqrt{p(1-p)q(1-q)}+pq)+n^2r(r+\varrho_2(1-r)).
\end{align*}
Combined, we see that
the difference in the objective function for $P$ as compared to $I_{2n}$ is}
\tb{\begin{align*}
\frac{1}{2}\e[\text{tr}(APBP^T)-\text{tr}(AB)]\geq&\binom{n}{2}(p^2+q^2)+(n^2-n)r^2\\
&-\binom{n}{2}(2\varrho_1\sqrt{p(1-p)q(1-q)}+2pq)-n^2r(r+\varrho_2(1-r))\\
=&\binom{n}{2}\left((p-q)^2-2\varrho_1\sqrt{p(1-p)q(1-q)}\right)-n^2\varrho_2r(1-r)-nr^2\\
=&n^2\left(\frac{1}{2}[(p-q)^2-2\varrho_1\sqrt{p(1-p)q(1-q)}]-\varrho_2r(1-r) \right)+O(n).
\end{align*} }
Numerous choices of the parameters in this model (for example, $p=0.8$, $q=0.2$, $r=0.1$, $\varrho_1=0.25$, $\varrho_2=0.3$) yields that 
$\e[\text{tr}(APBP^T)-\text{tr}(AB)]\geq cn^2+O(n)$ for a positive constant $c>0$ (in the example, $c=0.113$).
As $\text{tr}(APBP^T)-\text{tr}(AB)$ is highly concentrated about its expectation (see Appendix \ref{AP:TH5}), there is high probability that 
$\text{tr}(APBP^T)-\text{tr}(AB)>0$, and
$$\{I_{2n}\}\notin\argmax_{P\in \Pi(n)} \text{tr}(APBP^T)=\argmin_{P\in \Pi(n)} \|AP-PB\|_F;$$
i.e., $A$ and $B$ would not be $\gd$-matchable.
\end{example}

\subsection{Centering to recover matching}
\label{sec:cent}

\tb{In the previous example, we see that $\e(A)\neq\e(B)$ can effectively make $A$ and $B$ not $\gd$-matchable.
One way to recover the latent alignment in this heterogeneous setting is 
to transform the problem back into the homogeneous case, and} rather than matching $A$ and $B$, we would match $\tA=A-\e(A)$ and $\tB=B-\e(B)$; yielding once again $\e\tB=\e\tA$.
\tb{As the next theorem demonstrates, this is sufficient to a.s. recover $\td$-matchability under mild model conditions.
Before stating the theorem, we first must define some additional notation.
For $(A,B)\sim\HER$, and permutation $P\in\Pi(n)$, define the matrix $\ck_{A,B,P}$ via
\begin{align}
\label{eq:K}
\big[\mathcal{K}_{A,B.P}\big](u,v)=\text{Cov}\left[A(u,v),(PBP^T)(u,v)\right].
\end{align}
For each $P$ in $\Pi(n)$, define 
\begin{align}
\label{eq:xp}
\mathcal{X}_P:=\frac{1}{2}\left[\text{tr}\left(\ck_{A,B,I_n}\right)-\text{tr}\left(\ck_{A,B,P}\right)\right].
\end{align}}
\begin{theorem}
\label{thm:centeredmatched}
Let $(A,B)\sim \HER$ and consider $\tA=A-\e(A)$ and $\tB=B-\e(B)$.
For each $k\in[n]$, define
$$\Pi(n,k):=\left\{P\in\Pi(n)\text{ s.t. }\sum_{i=1}^n P(i,i)=n-k\right\}$$
If for all $k\in[n]$ and all $P\in \Pi(n,k)$, we have 
$$\frac{\mathcal{X}_P}{k}=\omega(\sqrt{n\log n}),$$ then
\begin{align*}
\p(A\text{ and }B\text{ are }\td\text{-matchable})=\p(\mathrm{argmin}_{P\in\Pi(n)}\|\tA-P\tB P^T\|_F=\{I_n\})\geq 1-e^{ -\omega(\log n)}.
\end{align*}
\end{theorem}
\noindent The proof of Theorem \ref{thm:centeredmatched} relies on a now standard application of McDiarmid's inequality and is similar to the proofs of analogous matchability results in \cite{JMLR:v15:lyzinski14a,rel,lyzinski2016information}; details of the proof can be found in Appendix \ref{AP:TH5}.

\begin{remark}
\label{remark:growth}
The growth condition on $\mathcal{X}_P$ in Theorem \ref{thm:centeredmatched}, namely
$\mathcal{X}_P/k=\omega(\sqrt{n\log n})$, 
is attempting to capture the necessary degree to which the entry-wise covariance matrix $\mathcal{C}$ needs to  be asymmetric.
If we define $\epsilon=\min_{u,v}\big[\mathcal{C}(A,B)\big](u,v)$, then from Eqs.\@ (\ref{eq:eps1}) and (\ref{eq:epnd}) we that for $P\in \Pi(n,k)$, 
\begin{align*}
\mathcal{X}_P\geq \frac{1}{2}\epsilon k\left(n-1-\frac{k}{2}\right),
\end{align*}
and if $\epsilon=\omega(\sqrt{\log n/n})$ then 
$\mathcal{X}_P/k=\omega(\sqrt{n\log n})$.
Constraining $\mathcal{X}_P$ globally and not entry-wise allows for more flexibility in applying the theorem to settings where some of the edges are very sparse or weakly correlated.
\end{remark}

\tb{Consider the growth condition on $\mathcal{X}_P$, namely
$\varrho p=\omega(\sqrt{\log n/n})$, in the $Q_1=Q_2=pJ_n$, $R=\varrho J_n$ homogeneous ER setting (wlog, assume $p<1/2$).
In this setting, as $\e(A)=\e(B)$, 
$A$ and $B$ are $\td$-matchable 
 iff $A$ and $B$ are $\gd$-matchable.
In the sparse setting of \cite[Theorem 1]{het2}, $\gd$-matchability is achieved with high probability when all of the following hold 
\begin{align}
\label{eq:p1}
p(1-p)(1-\varrho)&=O(1/\log n)\\
\label{eq:p2}
p(p+\varrho(1-p))&=O(1/\log n)\\
\label{eq:p3}
\frac{p(1-p)(1-\varrho)^2}{(p+\varrho(1-p))(1-p+p\varrho)}&=O(1/\log^3 n)\\
\label{eq:p4}
p(p+\varrho(1-p))&\geq \frac{\log n+\omega(1)}{n}.
\end{align}
Note that these conditions cannot simultaneously hold when
$$\frac{p}{\varrho}=\Omega(1).$$ 
Indeed if $p/\varrho=\Omega(1)$, Eq.\@ (\ref{eq:p2}) implies $p=o(1)$, and hence $\varrho=o(1)$. 
Therefore 
$$\frac{p(1-p)(1-\varrho)^2}{(p+\varrho(1-p))(1-p+p\varrho)}=\Omega(1),$$
contradicting Eq.\@ (\ref{eq:p3}).
In the $p/\varrho=\omega(1)$ setting, modulo the sparsity conditions, there is a $\gd$-matchability phase transition at $p^2+p\varrho= \frac{\log n}{n}$, and the corresponding rate achieved in Theorem \ref{thm:centeredmatched}, namely $\varrho p=\omega(\sqrt{\log n/n})$, is above this phase transition threshold.
We view this as the price paid in the Theorem for being able to handle both heterogeneous and homogeneous ER settings.}

\tb{In this setting, the growth condition of Theorem \ref{thm:centeredmatched}, $\varrho p=\omega(\sqrt{\log n/n})$, can hold in the dense setting $(p=\Theta(1))$ as well as when $p/\varrho=\Omega(1).$
In the dense setting of $p=\Theta(1)$, $\gd$-matchability transitions at $\varrho^2= \Theta(\frac{\log n}{n})$ \cite{lyzinski2016information,JMLR:v15:lyzinski14a}, which our Theorem recovers (asymptotically).}

\tb{In \cite{het2}, the authors establish a phase transition at $p^2+p\varrho= \frac{\log n}{n}$, providing the corresponding converse result that ensures no $\delta$-matchability if 
$p^2+p\varrho\leq \frac{\log n-\omega(1)}{n}$
for all dissimilarities $\delta$.
We do not derive a corresponding converse result to Theorem \ref{thm:centeredmatched} herein
(namely, a condition on $\mathcal{X}_P$ that ensures $A$ and $B$ are not $\delta$-matchable for any suitable $\delta$), as we are focused on how to practically recover $\delta$-matchability in the non-identically distributed setting; see Theorem \ref{thm:estimatedcenteredmatched}.  }

\subsection{Approximate centering to almost recover matching}
\label{sec:appcent}

Unfortunately, centering $A$ and $B$ by the true edge probability matrices $\e(A)$ and $\e(B)$ is impractical, as these model parameters are unknown in practice.
Our solution to this hurdle is to estimate the unknown $\e(A)$ and $\e(B)$ via Universal Singular Value Thresholding (USVT) \cite{chatterjee2015matrix}, and then \emph{approximately} center the networks via these estimates.

Our method for estimating $\e(A)$ and $\e(B)$ is based on the Universal Singular Value Thresholding (USVT) method of \cite{chatterjee2015matrix}, and
USVT applied in the present setting is outlined in Algorithm \ref{alg:USVT}.
\begin{algorithm}[t!]
  \begin{algorithmic}
    \STATE \textbf{Input}: Adjacency matrix $A$, threshold $t>0$;
 
\STATE{\bf 1.  } Let $A=\sum_{i=1}^n\sigma_iu_iv_i^T$ be the singular value decomposition of $A$, with singular values ordered via $\sigma_1\geq \sigma_2\geq \cdots\geq \sigma_n$;
\STATE{\bf 2.  } Let $T=\{\sigma_i>t\}$ be the set of singular values greater than the threshold $t$;
\tb{\STATE{\bf 3.  } Define $\breve{A}=\sum_{i=1}^{|T|}\sigma_iu_iv_i^T$}
\tb{\STATE\textbf{Output}: $\widehat{Q}\in[0,1]^{n\times n}$ defined via $\widehat{Q}(i,i)=0$ for all $i$, and for $i\neq j$
$$\widehat{Q}(i,j)=\begin{cases}
\breve{A}(i,j)&\text{ if }\breve{A}(i,j)\in[0,1]\\
0&\text{ if }\breve{A}(i,j) < 0\\
1&\text{ if }\breve{A}(i,j) > 1
\end{cases}$$}
\end{algorithmic}
\caption{USVT for estimating $Q$}
\label{alg:USVT}
\end{algorithm} 
If we estimate $\e(A)$ via $\hQ_1$ and $\e(B)$ via $\hQ_2$ using USVT, 
can we recover \tb{$\delta$-matchability} using the approximately centered matrices, $\hA=A-\hQ_1$ and $\hB=B-\hQ_2$, \tb{for a suitable $\delta$?}
Given the error introduced in estimating the edge probability matrices, the answer is unsurprisingly no, at least for the proof techniques we employ herein.
However, if we slightly weaken Definition \ref{def:matchable} to allow for a vanishing fraction of unmatched nodes, then we can recover an analogous result to Theorem \ref{thm:centeredmatched}.
This motivates the following definition:
\begin{definition}
\label{def:fmatchable}
Let $\mathbf{\delta}(\cdot,\cdot):\mathbb{R}^{n\times n}\times \mathbb{R}^{n\times n}\mapsto \mathbb{R}$ be a dissimilarity.  Consider random graphs $(A,B)\sim \HER$.
We say that $A$ and $B$ are $f_n,\delta$-matchable if 
$$\left\{\argmin_{P\in \Pi(n)} \delta(A, PBP^T)\right\}\cap\left\{ \cup_{i=f_n+1}^{\infty} \Pi(n,i)\right\}=\emptyset.$$
\end{definition}
\noindent Unwrapping Definition \ref{def:fmatchable}, we see that $f_n,\delta$-matchability is equivalent to any optimal permutation (under dissimilarity $\delta$) correctly recovering the labels of at least $n-f_n$ vertices across $A$ and $B$.
As the following theorem indicates, under mild model assumptions, $A$ and $B$ are \tb{$\Theta(\sqrt{n}),\hd$-matchable asymptotically almost surely, where the estimates $\widehat Q_A$ and $\widehat Q_B$ in Eq.\@ \ref{eq:delhat} are the USVT estimates. }
The proof of Theorem \ref{thm:estimatedcenteredmatched} can be found in Appendix \ref{pf:estimatedcenteredmatched}.
\begin{theorem}
\label{thm:estimatedcenteredmatched}
Let $(A,B)\sim \HER$ and further assume that for each $i=1,2$,
\begin{itemize}
\item[i.] There exists $r_i\in[0,1]$ such that $Q_i\leq r_i$ entry-wise.  Note that for each $n$, $r_i$ is fixed, though we allow $r_i$ to vary in $n$.
\item[ii.] We have that $\min_j \sum_\ell Q_i(j,\ell)=\omega(\log n)$.
\item[iii.] $Q_i$ is approximately low rank in that there exists a $d_i=o(n)$ such that $\sum_{j>d_i}s_j^2(Q_i)=O(nr_id_i)$ where $s_1(Q_i)\geq s_2(Q_i)\geq\cdots\geq s_n(Q_i)$ are the singular values of $Q_i$.
\end{itemize}
If for all $k\in[n]$ and all $P\in \Pi(n,k)$, we have that there exists $\delta=\delta(n)\geq 0$ such that
$$\frac{\mathcal{X}_P}{k}=\omega\left(n^{1/2+\delta}\sqrt{\log n}\right),$$ 
and for each $i=1,2$, there exists constants $a_i>0$ such that if $\hQ_i$ is the USVT estimate of $Q_i$ with threshold level $t_i=a_i\sqrt{nr_i}$ and if 
$$k_0=\omega\left(
\frac{n^{1-\delta}}{\sqrt{\log n}}(\sqrt{r_1}+\sqrt{r_2})\left(\sqrt{r_1 d_1}+\sqrt{r_2 d_2}\right) \right)$$
then $\hA=A-\hQ_1$ and $\hB=B-\hQ_2$ satisfy
\begin{align*}
\p&\left(
\mathrm{argmin}_{P\in\Pi(n)}\|\hA-P\hB P^T\|_F \bigcap 
\left\{ 
\cup_{i=k_0+1}^{\infty} \Pi(n,i)
\right\}
=\emptyset\right) \\
&=\p(A\text{ and }B\text{ are }k_0,\hd\text{-matchable})\geq 1-5n^{-2}.
\end{align*}
\end{theorem}

Let us take a moment to explore the assumptions in Theorem \ref{thm:estimatedcenteredmatched}.
Assumption $i.$ and $ii.$ control the allowable sparsity of the networks, ensuring that the minimum expected degree grows asymptotically faster than $\log n$.
\tb{If the mean expected degree was $o(\log n)$, then the graphs would be a.s.\@ disconnected \cite{boll}, and our proof techniques fail as $A$ and $B$ would no longer concentrate about $Q_1$ and $Q_2$ with high probability \cite{le2017concentration}.}
The rank assumption in $iii.$ is needed to control the accuracy of the USVT estimates of the unknown $Q_i$'s.  
Practically, smaller $d_i$ allow us to use suitable low-rank estimates $\widehat Q_i$ of $Q_i$ that are computationally easier to implement;
this is indeed the case in many common random graph models such as the Stochastic blockmodel \cite{sbm} (where often $d_i=O(1)$), random dot product graphs \cite{young2007random} and latent position random graphs \cite{hoff2002latent} (where $d_i$ is often taken to be $O(\mathrm{polylog}(n))$ \cite{udell2017nice}), among others.

\tb{If $R$ is bounded away from $0$ entry-wise, and each entry of $Q_i$ is $\Theta(r_i)$ (which is indeed the case in the oft adopted setting where each $Q_i=r_iM_i$ for a matrix $M_i$ with entries of order $\Theta(1)$) then $\epsilon$ as defined in Remark \ref{remark:growth} satisfies
$\epsilon=\Theta(\sqrt{r_1r_2})$.
We then have
$\mathcal{X}_P/k=\Omega(\sqrt{r_1r_2}\,n).$
From Eq. \ref{eq:k0} in the proof of Theorem \ref{thm:estimatedcenteredmatched}, we see that $k_0,\hd$-matchability is achieved here for 
$$k_0=\omega\left(
\frac{n^{1/2} }{\sqrt{r_1r_2}}(\sqrt{r_1}+\sqrt{r_2})\left(\sqrt{r_1 d_1}+\sqrt{r_2 d_2}\right) \right).$$
If, in addition $r_1/r_2=\Theta(1)$ and each $d_i=O(\text{polylog}(n))$, then up to a logarithmic factor $A$ and $B$ are $\Theta(\sqrt{n}),\hd$-matchable, and an oracle graph matching algorithm would properly align all but potentially a vanishing fraction of the nodes across the graphs.}

\begin{remark}
\label{rem:USVT}
	In Theorem \ref{thm:estimatedcenteredmatched}, we estimate $Q_i$ via $\hQ_i$ using USVT with threshold $t\approx\sqrt{n}$.
	In application, often suitable estimates of $Q_i$ can be obtained with rank $d_i$ of order $\log n$ or $\Theta(1)$ \cite{udell2017nice}, especially in the setting of latent space graph models.
	For the purposes of our proof approach, suitably good means that \tb{$\|\e(A)-\widehat Q_1\|_F=O(\sqrt{nr_id_i})$ (similarly for $B$).}
	We do not explore this model selection question further here (i.e., estimating a suitable rank rather than a threshold for our USVT estimates),
	as in applications often only a relatively small number of singular values are above the USVT threshold.
\end{remark}

\subsection{When to center?}
\label{sec:when?}

\tb{We have seen above that in the setting where both $Q_1\neq Q_2$ and $A,B$ are not $\gd$-matchable, centering $A$ and $B$ via $\tA$ and $\tB$ can recover $\td$-matchability by ameliorating the effect of the differing $Q_i$'s.
Moreover, approximately centering by $\hQ_1$ and $\hQ_2$ theoretically recovers $\hd$-matchability for all but a vanishing fraction of the vertices.
A natural question is in the case when $Q_1=Q_2$, does Theorem \ref{thm:estimatedcenteredmatched} imply that  a.s.\@ perfect $\gd$-matchability is potentially lost when USVT centering is performed unnecessarily?}

\begin{figure}[t!]
	\centering
	\includegraphics[width=0.6\textwidth]{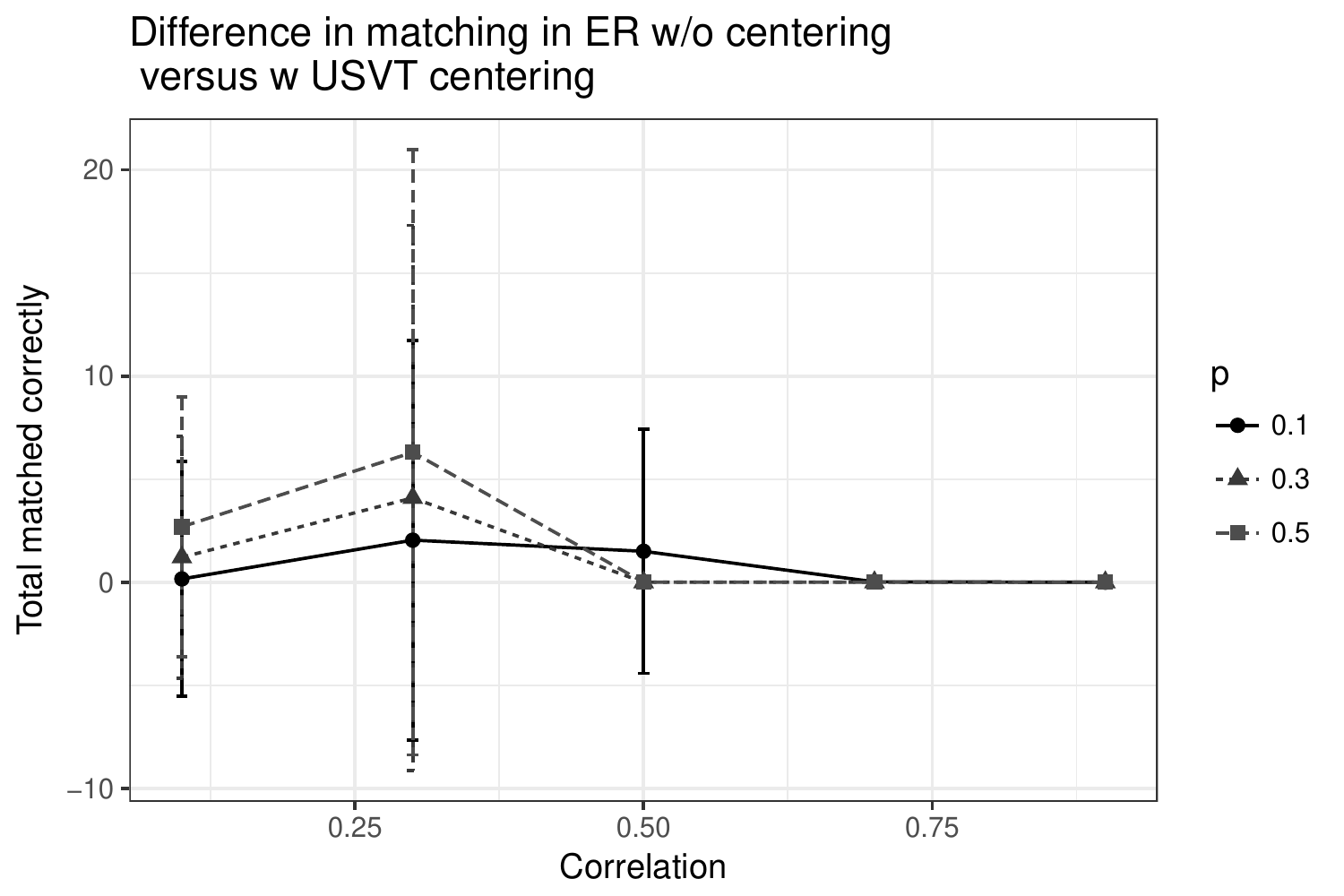} 
	\caption{
	\tb{We plot the mean ($\pm$ 1 s.d.) of 
\texttt{FAQ}$(A,B,I_n)-$\texttt{FAQ}$(A-\hQ_1,B-\hQ_2,I_n)$ over $p\in\{0.1,0.3,0.5\}$
and $\varrho\in\{0.1,0.3,0.5,0.7,0.9\}$.
	The parameters considered are $p\in\{0.1,0.3,0.5\}$
	and $\varrho\in\{0.1,0.3,0.5,0.7,0.9\}$.}}
	\label{fig:ercent}
\end{figure}

\tb{Consider the following simple example, where $Q_1=Q_2=pJ_n$ and $R=\varrho J_n$, $n=100$, and we vary $p\in\{0.1,0.3,0.5\}$
and $\varrho\in\{0.1,0.3,0.5,0.7,0.9\}$.
In this example, there is no need to center before matching, and the variability introduced by estimating the $Q_i$'s could potentially cause 
$$
\mathrm{argmin}_{P\in\Pi(n)}\|\hA-P\hB P^T\|_F\neq\{I_n\}.
$$
Fortunately, at least in this example we see this is not the case (see Figure \ref{fig:ercent}).
As matching these graphs exactly (i.e., finding the argmin of the GMP) is computationally challenging, we use as a surrogate for $\gd$-matchability 
and 
$\hd$-matchability
whether the true alignment is a local (rather than global) minima of the GMP before and after centering.
To test this, we match the graph pairs (USVT centered and uncentered) using the constrained gradient-ascent based graph matching algorithm, \texttt{FAQ} \cite{FAQ}, initialized at the true correspondence $I_n$.
While \texttt{FAQ} is not guaranteed to terminate at a local minima, if it terminate at $I_n$, then that is evidence in support of $I_n$'s local optimality.
Moreover, if \texttt{FAQ} does not terminate at $I_n$, then $I_n$ is not a local minima of the GMP.}

Letting \texttt{FAQ}$(A,B,I_n)$ (resp., \texttt{FAQ}$(A-\hQ_1,B-\hQ_2,I_n)$) denote the number of vertices correctly matched by \texttt{FAQ} initiated at $I_n$ when $A$ and $B$ are matched directly (resp., when $A$ and $B$ are USVT centered before they are matched).
In Figure \ref{fig:ercent}, we plot the mean ($\pm$ 1 s.d.) of 
\texttt{FAQ}$(A,B,I_n)-$\texttt{FAQ}$(A-\hQ_1,B-\hQ_2,I_n)$ over $p\in\{0.1,0.3,0.5\}$
and $\varrho\in\{0.1,0.3,0.5,0.7,0.9\}$.
From the figure, we see that there is no significant performance lost by centering the graphs first.
Indeed, in highly structured/low-rank settings (e.g., homogeneous ER or SBM), we can obtain high-fidelity estimates of the individual entries of the $Q_i$ versus the global estimates used in current proof.
These local estimates (which can also be obtained by non-spectral methods, e.g., in the $\mathrm{ER}(n,p)$ case we can use $\hat p=\sum_{i<j}A(i,j)/\binom{n}{2}$) should allow for significantly less error to be introduced in the estimation of the $Q_i$'s and will allow for little--to--no theoretic degradation due to centering.
As we are more focused on the general $Q_i$ case, we do not pursue this further here. 

\subsection{Core Matchability}
\label{sec:cores}

\tb{Often in applications, only a fraction of the vertices in $A$ possess a latent matched pair in $B$.
We will denote those vertices that have a latent match across graphs as the core vertices (denoted $\mathcal{C}$), 
and we will denote those vertices that do not have a latent match across graphs as the junk vertices (denoted $\mathcal{J}$).
In this section, we seek to further understand the ability of an oracle graph matching procedure to correctly match the cores across graphs.
This motivates the following definition of core-matchability.}
\tb{\begin{definition}
\label{def:corematch}
Let $(A,B)\sim\HER$, and consider a partition of the vertex sets into core and junk vertices,  
$$V=\mathcal{C}\cup \mathcal{J}.$$ 
Define $\mathcal{P}_n^\mathcal{C}:=\{P\in\Pi(n)\text{ s.t. } \forall i \in \mathcal{C}, P(i,i)=1\}$ to be the set of core matching permutations.
For dissimilarity $\delta$, we say that $A$ and $B$ are \emph{core} $\delta$-\emph{matchable} if
$$\argmin_{P\in\Pi(n)}\delta(A,PBP^T)\subset \mathcal{P}_n^\mathcal{C};$$
i.e., if any optimal permutation aligning $A$ and $B$ under $\delta$ perfectly matches the cores across networks. 
\end{definition}}

\tb{If we consider $(A,B)\sim\HER$ with $R$ of the form $R=R_c\oplus 0_{n_j}$ where $R_c\in(0,1]^{n_c\times n_c}$, then it is reasonable to define $\mathcal{C}=[n_c]$ and $\mathcal{J}=[n]\setminus[n_c]$.
Indeed, under these model assumptions $R(u,v)=0$ if either $u$ or $v$ is in $\mathcal{J}$.}

Completely analogously to the setting considered in Example \ref{ex:bad}, it is immediate that $A$ and $B$ need not be asymptotically almost surely core $\gd$-matchable even with non-vanishing core correlation in $R$.  
\tb{Indeed, as in Example \ref{ex:bad}, $Q_1$ and $Q_2$ can be chosen to effectively obfuscate the true alignment amongst the core vertices.}
Mimicking the results of Theorem \ref{thm:centeredmatched}, centering $\tA=A-\e(A)$ and $\tB=B-\e(B)$ again a.s.\@ recovers core $\td$-matchability of $A$ and $B$ under mild model assumptions.  
The proof of Theorem \ref{thm:matchedcorescentered} is contained in Appendix \ref{AP:stochmatchedcores}.
\tb{\begin{theorem}
\label{thm:matchedcorescentered}
Let $(A,B)\sim \HER$ and consider $\tA=A-\e(A)$ and $\tB=B-\e(B)$.
Suppose that $R$ is of the form $R=R_c\oplus 0_{n_j}$ where $R_c\in(0,1]^{n_c\times n_c}$, and for each $P\in\Pi(n)$, let $\mathcal{X}_P$ be defined as in Eq.\@ (\ref{eq:xp}).
If for all $Q\in\{P\in\Pi(n)\text{ s.t. }\sum_{i=1}^{n_c}P(i,i)=n_c-k\}$ we have that 
\begin{align}
\label{eq:core1}
\mathcal{X}_Q/k&=\omega\left( \sqrt{n_c\log{n_c}}\right),\\
\label{eq:core2}
\mathcal{X}_Q/\sqrt{k}&=\omega\left( \sqrt{n_c\,n_j\log{n_j}}\right),
\end{align}
and also if $n_c\geq n_j$, then
$$\Pr\left(A\text{ and }B\text{ are core }\td\text{-matchable}\right)\geq1- 2e^{-\omega(\log n_c)}.$$
\end{theorem}}

\tb{\noindent As before, if we define
$\epsilon:=\min_{u,v\in[n_c]} Cov(A(u,v),B(u,v)),$
then for $Q\in\Pi(n,k)$,
$\mathcal{X}_Q/k=\Omega(\epsilon n_c)$.
If, in addition,
\begin{align*}
\epsilon&=\omega\left( \sqrt{\frac{n_j\log(n_j)}{n_c}}\right)\\
\epsilon&=\omega\left( \sqrt{\frac{\log(n_c)}{n_c}}\right),
\end{align*}
then Eqs. (\ref{eq:core1}) and (\ref{eq:core2}) hold and core $\td$-matchability is recovered.
In the event that $\epsilon$ is $\Theta(1)$ or $\Theta(1/\text{polylog}(n)),$ then Theorem \ref{thm:matchedcorescentered} implies that $A$ and $B$ are core $\delta$-matchable in the presence of nearly linear junk, \emph{with arbitrary junk structure}.
This result extends and generalizes the results in \cite{kazemi2015can} to the non-homogeneous ER setting.}

As before, if the unknown $\e(A)$ and $\e(B)$ are estimated via USVT, then we recover partial core matchability.
Before formalizing this, we first need the following extension of Definition \ref{def:fmatchable} to the core-junk setting.
\begin{definition}
\label{def:corefmatchable}
Let $\delta:\mathbb{R}^{n\times n}\times \mathbb{R}^{n\times n}\mapsto\mathbb{R}$ be a dissimilarity.
Let $(A,B)\sim\HER$, and consider a partition of the vertex sets into core and junk vertices,  
$$V=\mathcal{C}\cup \mathcal{J}.$$ 
Define
$$\mathcal{P}^{\mathcal{C}}_{n,k}:=\{P\in\Pi(n)\text{ s.t. }\sum_{i\in\mathcal{C}} P(i,i)=|\mathcal{C}|-k\}$$
We say that $A$ and $B$ are core $f_n,\delta$-matchable if 
$$\left\{\argmin_{P\in \Pi(n)} \delta(A,PBP^T)b\right\}\cap\left\{ \cup_{k=f_n+1}^{\infty} \mathcal{P}^{\mathcal{C}}_{n,k}\right\}=\emptyset.$$
\end{definition}

The following Theorem provides the analogue of Theorem \ref{thm:estimatedcenteredmatched} in the core-junk setting.  
Note that the proof is completely analogous to that in Theorem \ref{thm:estimatedcenteredmatched}, and so is omitted.
\begin{theorem}
\label{thm:coreestcenteredmatched}
Let $(A,B)\sim \HER$ and consider $\hA=A-\widehat Q_1$ and $\hB=B-\widehat Q_2$.
Suppose that $R$ is of the form $R=R_c\oplus 0_{n_j}$ where $R_c\in(0,1]^{n_c\times n_c}$.
With the assumptions on $Q_1$ and $Q_2$ from Theorem \ref{thm:estimatedcenteredmatched}, 
and assume 
for each $Q\in\{P\in\Pi(n)\text{ s.t. }\sum_{i=1}^{n_c}P(i,i)=n_c-k\}$, 
\begin{align}
\label{eq:core11}
\mathcal{X}_Q/k&=\omega\left( n_c^{1/2+\delta}\sqrt{\log{n_c}}\right),\\
\label{eq:core21}
\mathcal{X}_Q/\sqrt{k}&=\omega\left( \sqrt{n_c\,n_j\log{n_j}}\right),
\end{align}
and also assume $n_c\geq n_j$.
For each $i=1,2$, there exists constants $a_i>0$ such that if $\hQ_i$ is the USVT estimate of $Q_i$ with threshold level $t_i=a_i\sqrt{n\rho_i}$ and if 
$$k_0=\omega\left(
\frac{n_c^{1-\delta}}{\sqrt{\log n_c}}(\sqrt{\rho_1}+\sqrt{\rho_2})\left(\sqrt{\rho_1 d_1}+\sqrt{\rho_2 d_2}\right) \right),$$
then 
$$\p(A\text{ and }B\text{ are core }k_0,\hd\text{-matchable})\geq 1-5n^{-2}.$$
\end{theorem}

\section{Simulations and Experiments}
\label{sec:expsim}
In the following sections, we explore the impact on \emph{graph matchability} of USVT centering in both simulated and real data settings.
We note here that precisely determining the level of $\gd$, $\hd$ and $\td$-matchability is infeasible for even modestly sized networks, as this would require exactly solving the NP-hard graph matching problem.
To circumvent this, we instead match our networks using the \texttt{FAQ} algorithm of \cite{FAQ} initialized at a variety of starting points including the true correspondence.  
As it is a Frank-Wolfe \cite{FW} based algorithm, if \texttt{FAQ} terminates at the true correspondence (or at a permutation which matches a high percentage of the vertices), then the true correspondence is an \tb{estimated} local minima of the graph matching problem.
\tb{Moreover, if \texttt{FAQ} is initialized at the true correspondence and does not terminate at the true correspondence, then the true correspondence is not a local minima. }
Comparing objective function values across \tb{estimated} local minima then allows us to approximately gauge the global optimality of the true correspondence.
While this is not the same as directly finding the global minima desired in the definition of matchability, it nonetheless provides a useful, principled heuristic for empirically studying both matchability and deviations there from.

\begin{figure}[t!]
	\centering
	\includegraphics[width=0.7\textwidth]{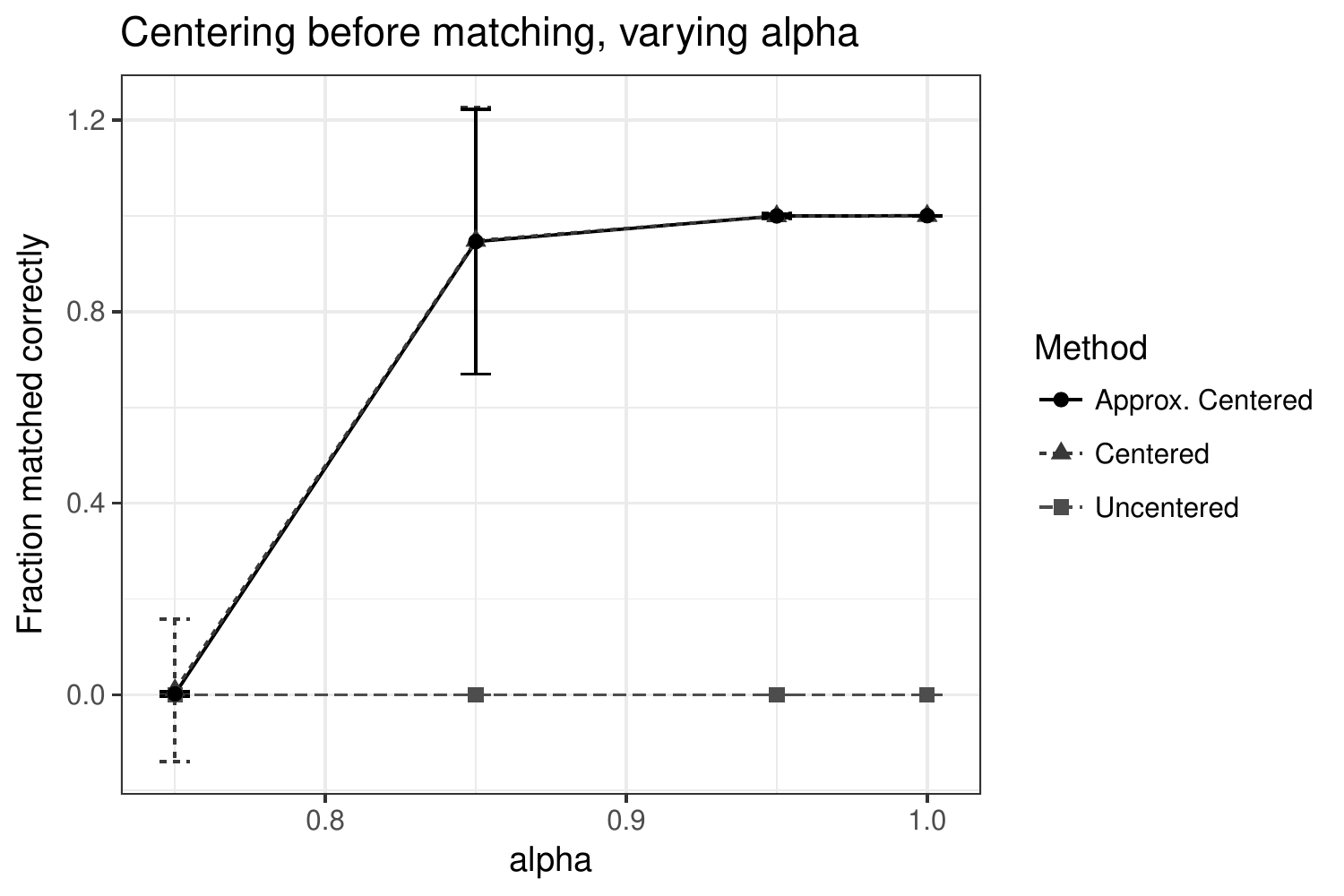}
	\caption{Fraction correctly matched by \texttt{FAQ} \tb{$\pm2$s.d.} (optimized over two different initializations: $I_{2n}$ and at $\tilde I_{2n}$) 
versus 
	$\alpha\in(0.75, 0.85, 0.95, 1)$ when matching i. $A$ and $B$ directly (labeled ``Uncentered''); ii. $\tA=A-Q_1$ and $\tB=B-Q_2$ (labeled ``Centered''); and $\hA=A-\hQ_1$ and $\hB=B-\hQ_2$ (labeled ``Approx. Centered'').
	Here, $\alpha$ captures the level of correlation between $A$ and $B$ (higher $\alpha$ means more correlation), $n=150$ is fixed and results are averaged over $100$ Monte Carlo replicates.}
	\label{fig:centered1}
\end{figure}
\subsection{Simulation}
\label{sec:censim}
\begin{figure}[t!]
	\centering
	\includegraphics[width=0.7\textwidth]{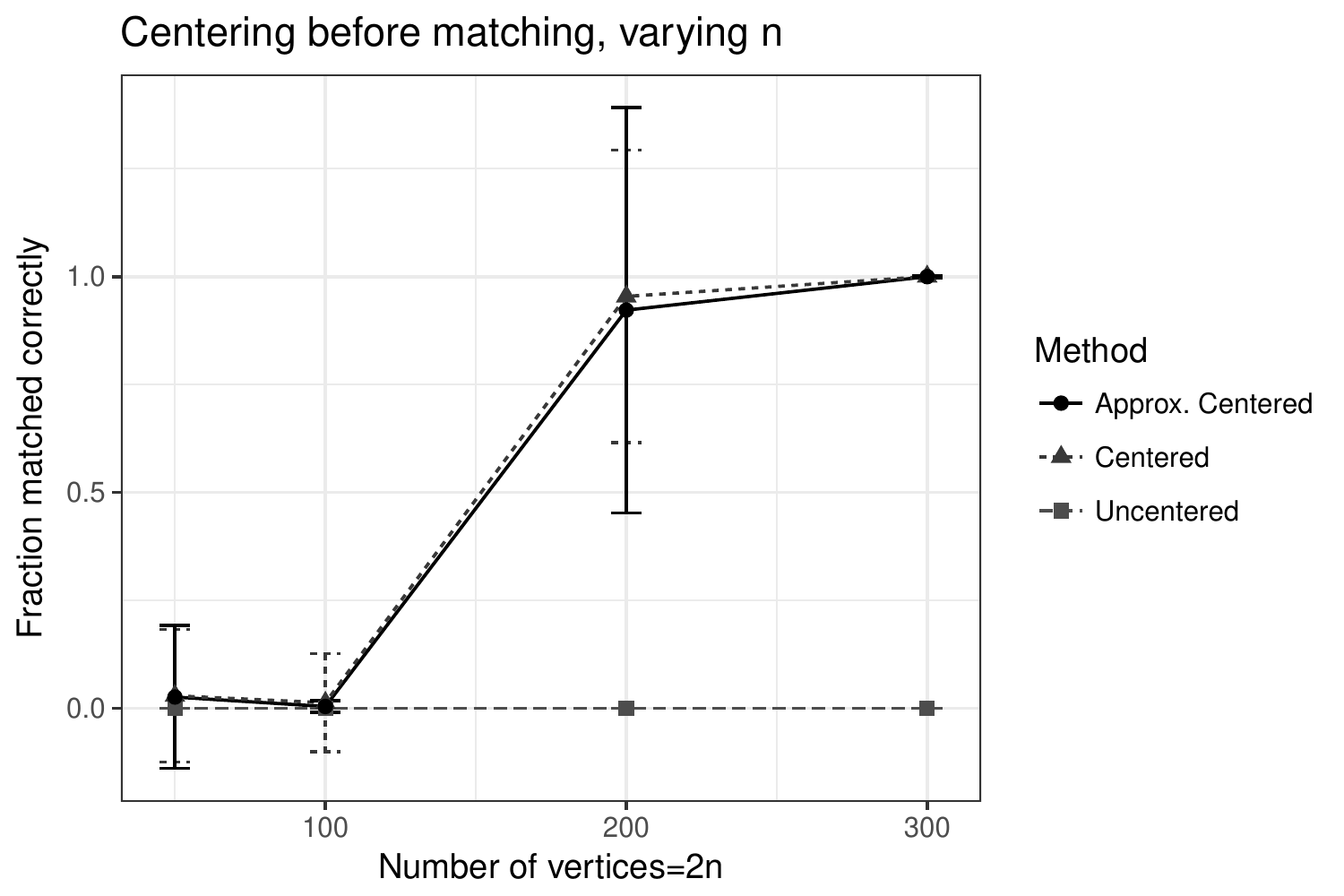}
	\caption{Fraction correctly matched by \texttt{FAQ} \tb{$\pm2$s.d.} (optimized over two different initializations---$I_{2n}$ and at $\tilde I_{2n}$) versus 
	$n\in(25,50,100,150)$ when matching i. $A$ and $B$ directly (labeled ``Uncentered''); ii. $\tA=A-Q_1$ and $\tB=B-Q_2$ (labeled ``Centered''); and $\hA=A-\hQ_1$ and $\hB=B-\hQ_2$ (labeled ``Approx. Centered'').
	Here, $\alpha=1$ is fixed and results are averaged over $100$ Monte Carlo replicates.}
	\label{fig:centered2}
\end{figure}

To explore the utility of USVT centering as a graph matching preprocessing step, we consider the following experiment.  We let $(A,B)\sim\HER$ with 
$$Q_1=\begin{pmatrix}
0.8J_n& 0.1J_n\\
0.1J_n& 0.2J_n
\end{pmatrix}
\text{ and }
Q_2=\begin{pmatrix}
0.2J_n& 0.1J_n\\
0.1J_n& 0.8J_n
\end{pmatrix}.$$
and $R$ of the form
$$R=\alpha\begin{pmatrix}
0.25J_n& 0.3J_n\\
0.3J_n& 0.25J_n
\end{pmatrix},$$
and we use the \texttt{FAQ} algorithm of \cite{FAQ} to match i. $A$ and $B$ directly (labeled ``Uncentered'' in Figures \ref{fig:centered1}--\ref{fig:centered2}); ii. $\tA=A-Q_1$ and $\tB=B-Q_2$ (labeled ``Centered'' in Figures \ref{fig:centered1}--\ref{fig:centered2}); and $\hA=A-\hQ_1$ and $\hB=B-\hQ_2$ (labeled ``Approx. Centered'' in Figures \ref{fig:centered1}--\ref{fig:centered2}).
In each figure, we initialize the \texttt{FAQ} algorithm at $I_{2n}$---i.e., at the true latent alignment---and at $\tilde I_{2n}:=\begin{pmatrix}0&I_{n}\\
I_{n}&0  \end{pmatrix}$---i.e., at the alignment completely confusing blocks one and two across networks.
We plot the mean fraction of vertices matched correctly \tb{($\pm 2$s.d.)} by \texttt{FAQ} at the starting point that achieves the lowest graph matching objective function score (averaged over $100$ Monte Carlo replicates).
As mentioned above, if the fraction matched correctly is less than $1$, then the true alignment is not a local minimum of the graph matching objective function and the graphs are not $\gd$, $\hd$, or $\td$-matchable (depending on what input \texttt{FAQ} is matching).

In Figure \ref{fig:centered1}, we consider $n=150$---i.e., the graphs are size $300$---and in the USVT estimates $\hQ_i$ we used $t=2.01\sqrt{300*0.8*0.2}$ as suggested in \cite{chatterjee2015matrix}.
We plot $\alpha\in(0.75, 0.85, 0.95, 1)$ versus the the mean fraction of vertices matched correctly \tb{($\pm 2$s.d.)}
It is unsurprising in light of Example \ref{ex:bad}, that $A$ and $B$ are not only directly not $\gd$-matchable, but that the alignment found by \texttt{FAQ} matches none of the vertices correctly across graphs.
Also of note in the figure is that as $\alpha$ increases, the oracle centered graphs $\tA$ and $\tB$ appear to be nearly $\td$-matchable (in that the estimated local minimum found by \texttt{FAQ} is close to $I_{300}$).
The steep performance drop off as $\alpha$ decreases is a consequence of the fact that in low-correlation regimes (low for a given $n$), $\td$-matchability is often not recovered even through centering.
Performance in the approximately centered case tracks performance in the centered case, with the USVT centering recovering the gains of the oracle centering.
This empirically suggests that the $k_0$ lower bound in Theorem \ref{thm:estimatedcenteredmatched} is not sharp, as USVT centering recovers full $\hd$-matchability in the high $\alpha$ regime.
\tb{We surmise that if $Q_i$ is truly low-rank, USVT centering and true centering will recover perfect $\hd$- and $\td$-matchability respectively as $n$ increases.}

In Figure \ref{fig:centered2} we repeat the above simulation with $\alpha=1$ fixed and $n\in(25,50,100,150)$. 
Using $t=2.01\sqrt{(2n*0.8*0.2)}$ as the USVT threshold, we again plot the mean fraction of vertices matched correctly ($\pm2$ s.d.) versus $n$.
As before, without centering the optimal alignment found by \texttt{FAQ} matches very few of the vertices correctly across graphs.
The performance increase in the oracle centered setting as $n$ increases is a consequence of the $\td$-matchability (in Theorem \ref{thm:centeredmatched}) being an asymptotic result; indeed, we should not expect correlated small networks to be almost surely $\td$-matchable even with the oracle centering.  
We note, however, that $n=150$ (i.e., graph order $=300$) here is sufficient for the asymptotically perfect $\td$-matchability to be recovered.
Again we see that performance in the approximately centered case tracks performance in the centered case, with the USVT centering achieving almost all of the gains of the oracle centering.


\subsection{Twitter data}
\label{sec:data}
\begin{figure}[t!]
	\centering
\begin{tabular}{ccc}
	\includegraphics[width=.2\textwidth]{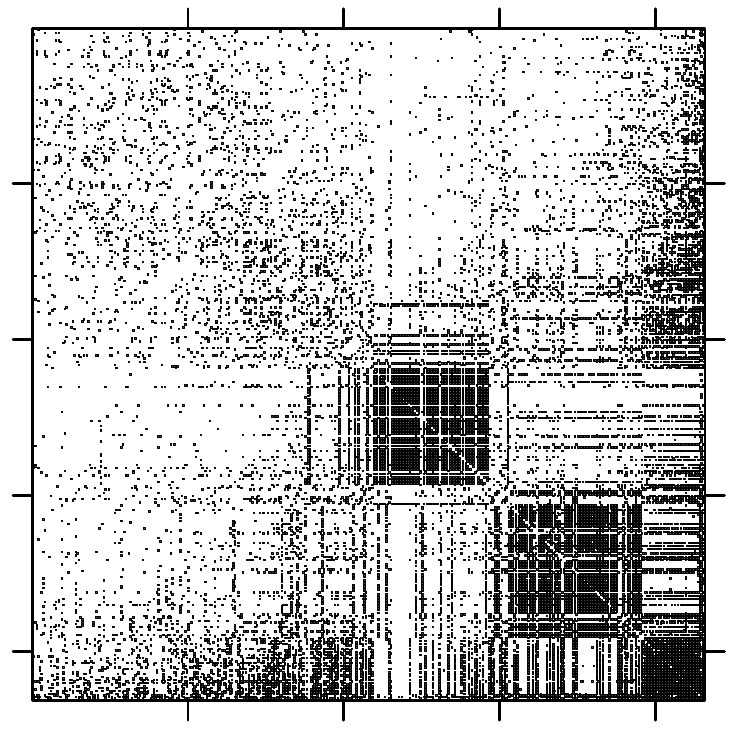} &
	\includegraphics[width=.2\textwidth]{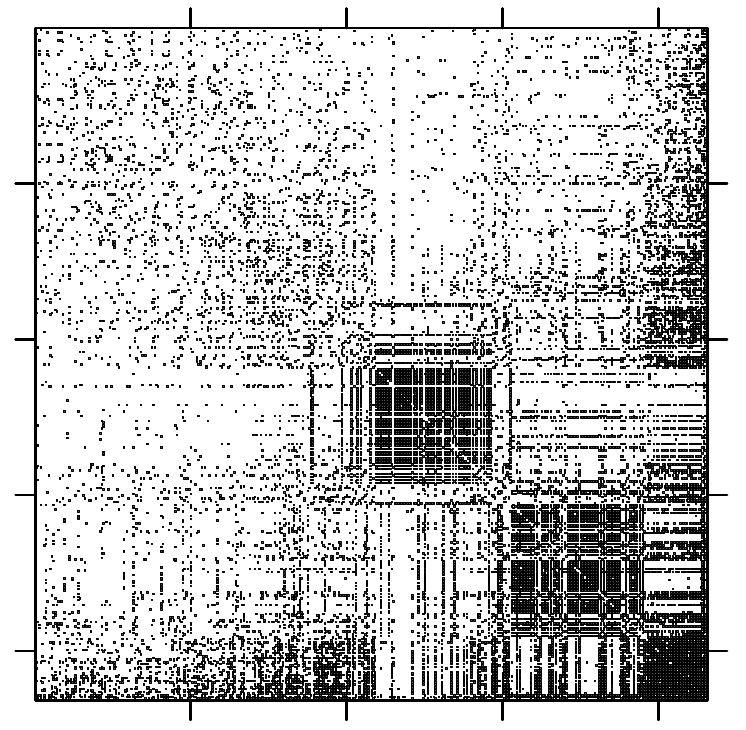}&
	\includegraphics[width=.5\textwidth]{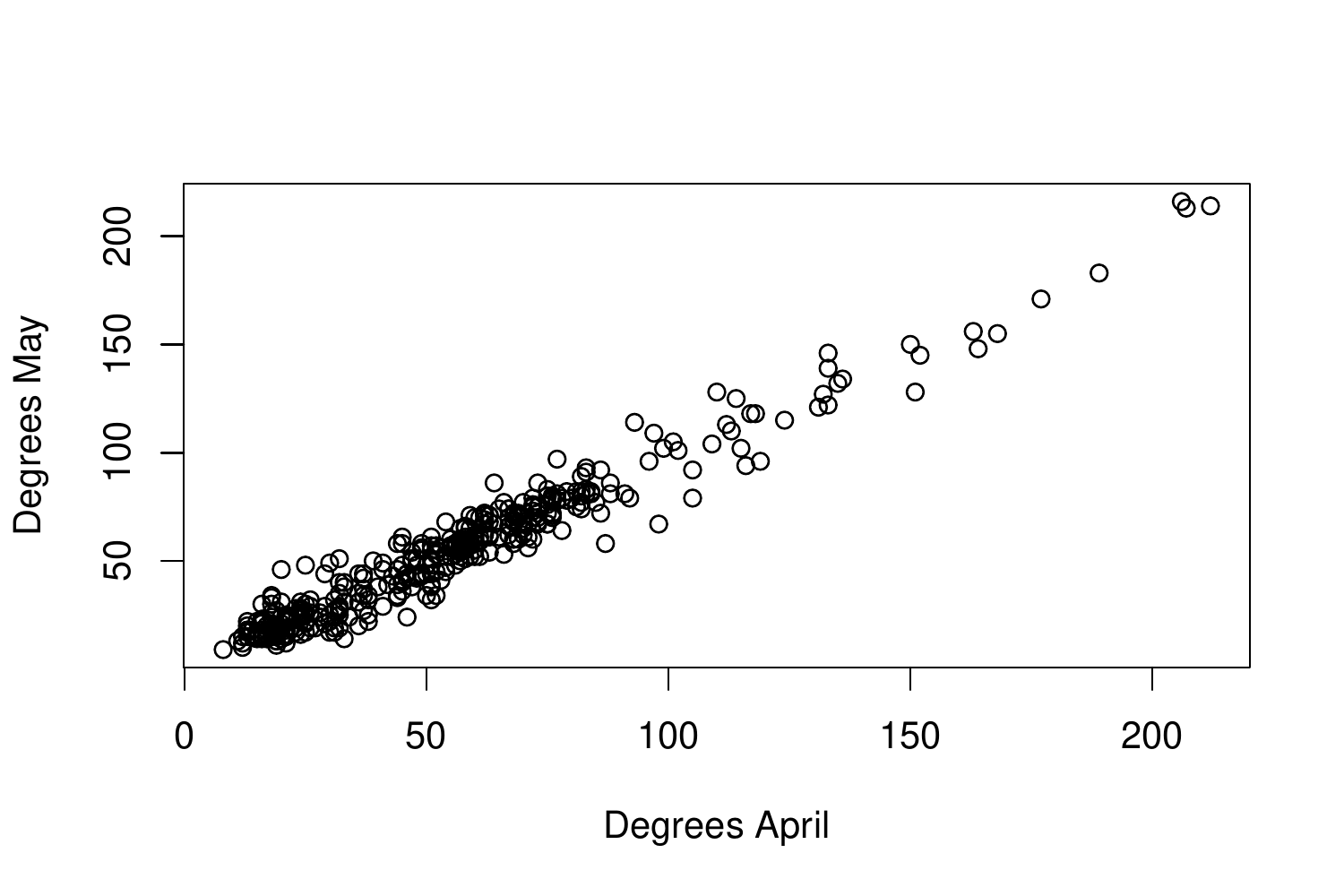} 
	 \\
	April & May & Degree Comparison
\end{tabular}
	\caption{In the left two panels, we plot the  adjacency matrices of aligned Twitter graphs from April and May.
	In the right panel, we plot the degrees of each vertex in April versus the degrees of the same vertex in May. 
	The vertices were both sorted according to ascending degree in the April graph.}
	\label{fig:twitsame}
\end{figure}

In order to analyze the impact of USVT centering with real data, we considered two graphs derived from Twitter.\footnote{These graphs were provided as part of the DARPA XDATA project.}
The two graphs are based on the most active twitter users from April and May 2014.
The graphs are \tb{unweighted} with an edge between users if a user mentioned another user during the given month.
After keeping only the largest common connected component, the number of users was 431 in each graph.

As can be seen in Figure~\ref{fig:twitsame}, many vertices in this data set have very similar connectivity patterns.
Indeed, the empirical Pearson correlation between the entries in the two adjacency matrices is \tb{$\approx 0.72$.}
It is not surprising that the similar graph topology across networks leads to good performance when matching the graphs without centering.
Repeating the experiment in Section \ref{sec:censim}---i.e., matching the adjacency matrices $A$ (``April'') and $B$ (``May'') using \texttt{FAQ} initialized at $I_{431}$---yields \tb{$323$ }vertices correctly aligned across networks.  
Although the true correspondence is not optimal (according to the GM objective function), the estimated local optimal correspondence does match $\approx 75\%$ of the vertices correctly across networks without the need for centering.
It is worth noting that preprocessing the data via USVT centering before matching again yields \tb{$323$} vertices correctly matched by \texttt{FAQ} initialized at $I_{431}$ (with a suitable USVT threshold here being \tb{$2\sqrt{431}$).}
This suggests that the centering procedure does not hurt performance when the graph topologies are similar across networks, and as we will demonstrate below, can significantly increase performance when the graph topologies differ across networks.
\tb{We do note here that we do not zero out the diagonal of $\hQ_i$ in the USVT step in this real data example, as here, hollow $\hQ_i$'s led to significantly worse performance than the non-hollow $\hQ_i$'s.}

\begin{figure}[t!]
	
\begin{tabular}{cc}
	\includegraphics[width=.5\textwidth]{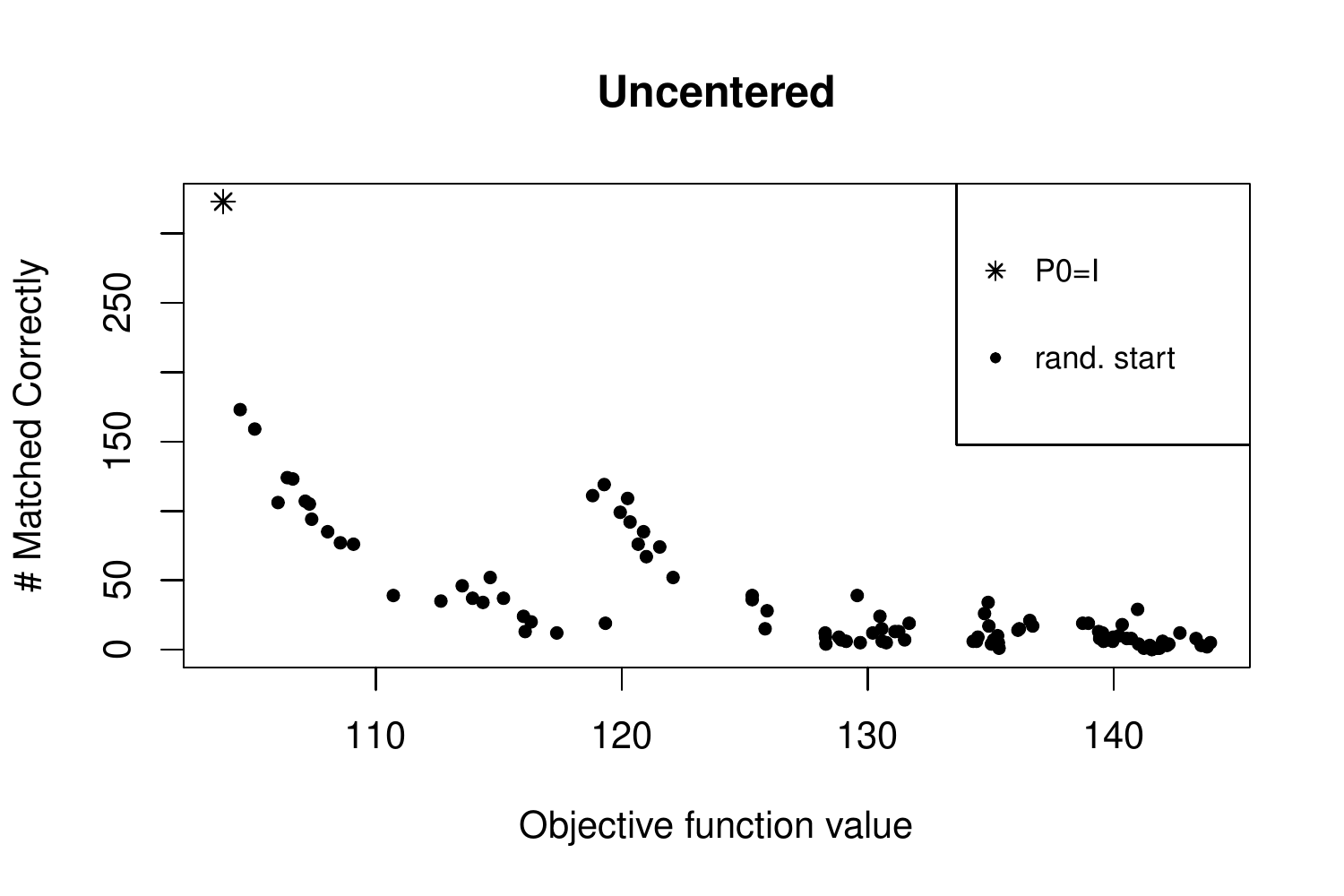} &
	\includegraphics[width=.5\textwidth]{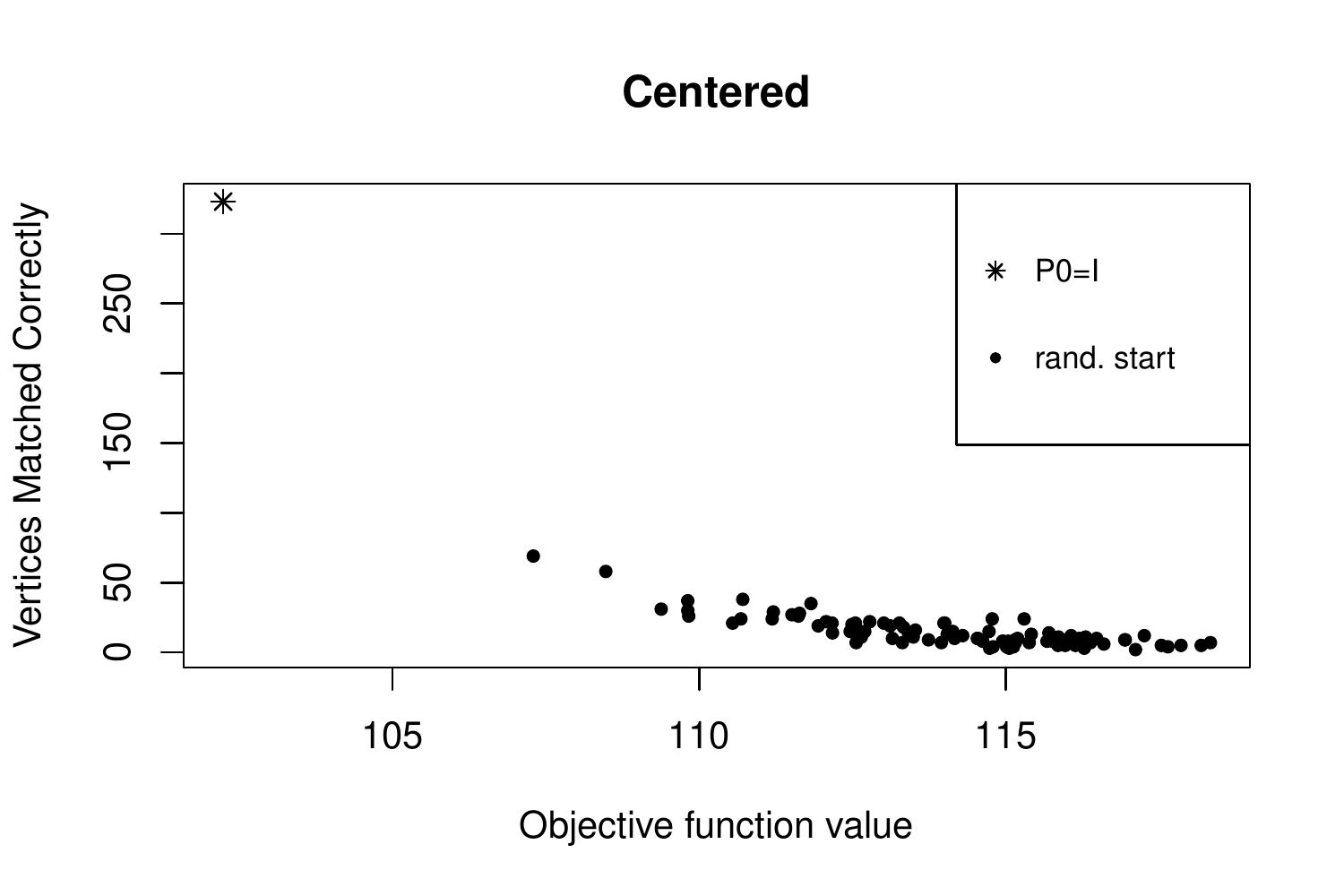}
\end{tabular}
	\caption{We plot the number of vertices in the April-May Twitter graphs correctly matched by \texttt{FAQ} versus the graph matching objective function value.
	In each panel (on the left matching $A$ and $B$, and on the right $\hA$ and $\hB$), we initialize \texttt{FAQ} at 100 different starting points:  once at $I_n$ (labeled ``P0=I'' in the legend) and the rest at random restarts (labeled ``rand. start'' in the legend).}
	\label{fig:twitcool}
\end{figure}

While both the centered and uncentered graphs are highly $\hd$- and $\gd$-matchable respectively, centering does have a very interesting algorithmic effect here; see Figure \ref{fig:twitcool}.
In the figure, we plot the number of vertices in the April-May Twitter graphs correctly matched by \texttt{FAQ} versus the graph matching objective function value.
In each panel (on the left matching $A$ and $B$, and on the right $\hA$ and $\hB$), we initialize \texttt{FAQ} at 100 different starting points:  once at $I_n$ (labeled ``P0=I'' in the legend) and the rest at random permutation restarts (labeled ``rand. start'' in the legend).
The figure suggests that centering has the effect of creating a more stable objective function gap between the estimated optimal permutation and suboptimal alternatives.
In a setting where multiple random restarts are possible---and needed---to recover an unknown latent alignment, this suggests that the optimal alignment is perhaps more easily recognized in the centered graph regime, and hence online stopping criterion more easily implementable.

\begin{figure}[t!]
	\centering
\begin{tabular}{cc}
	\includegraphics[width=.5\textwidth]{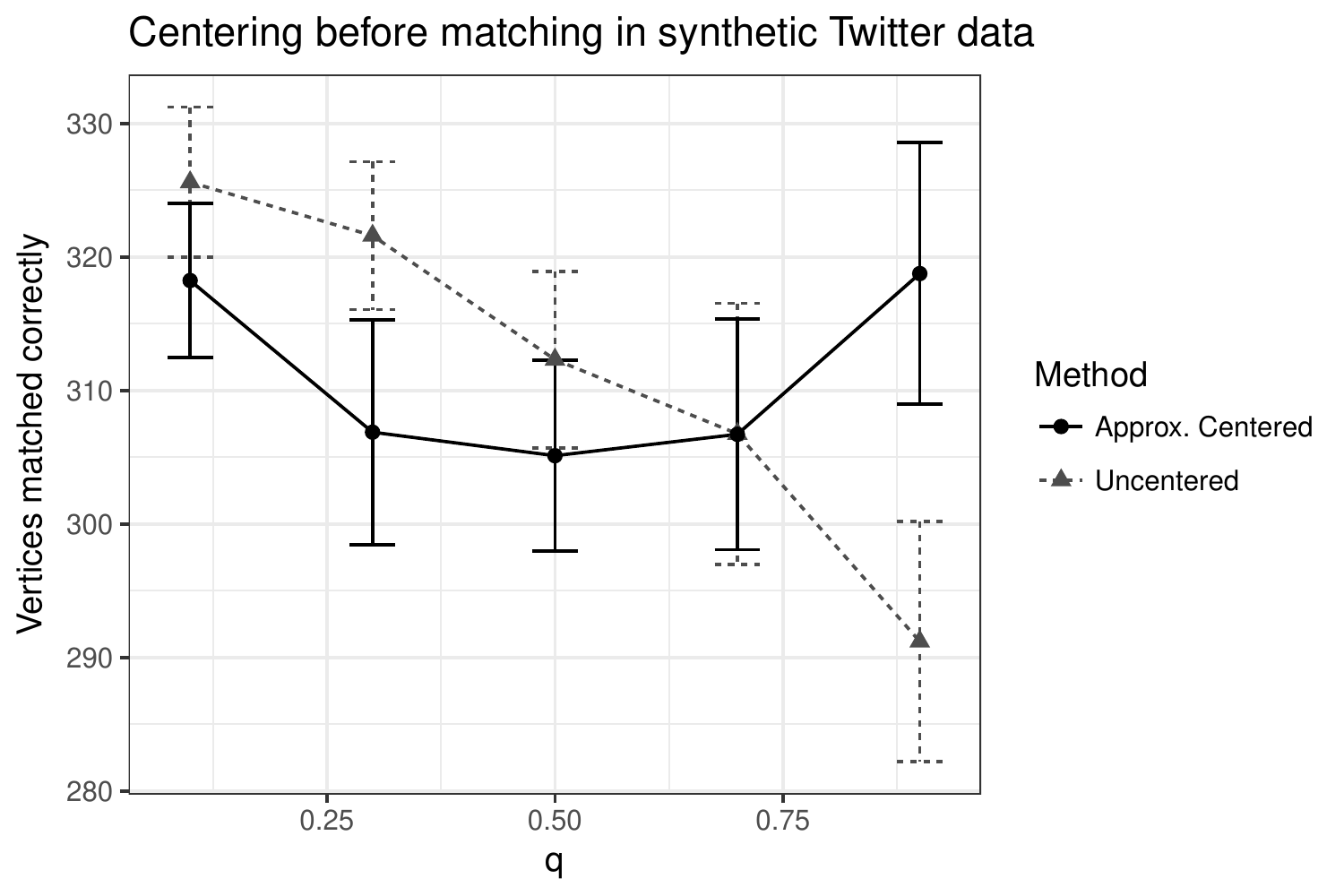} &
	\includegraphics[width=.5\textwidth]{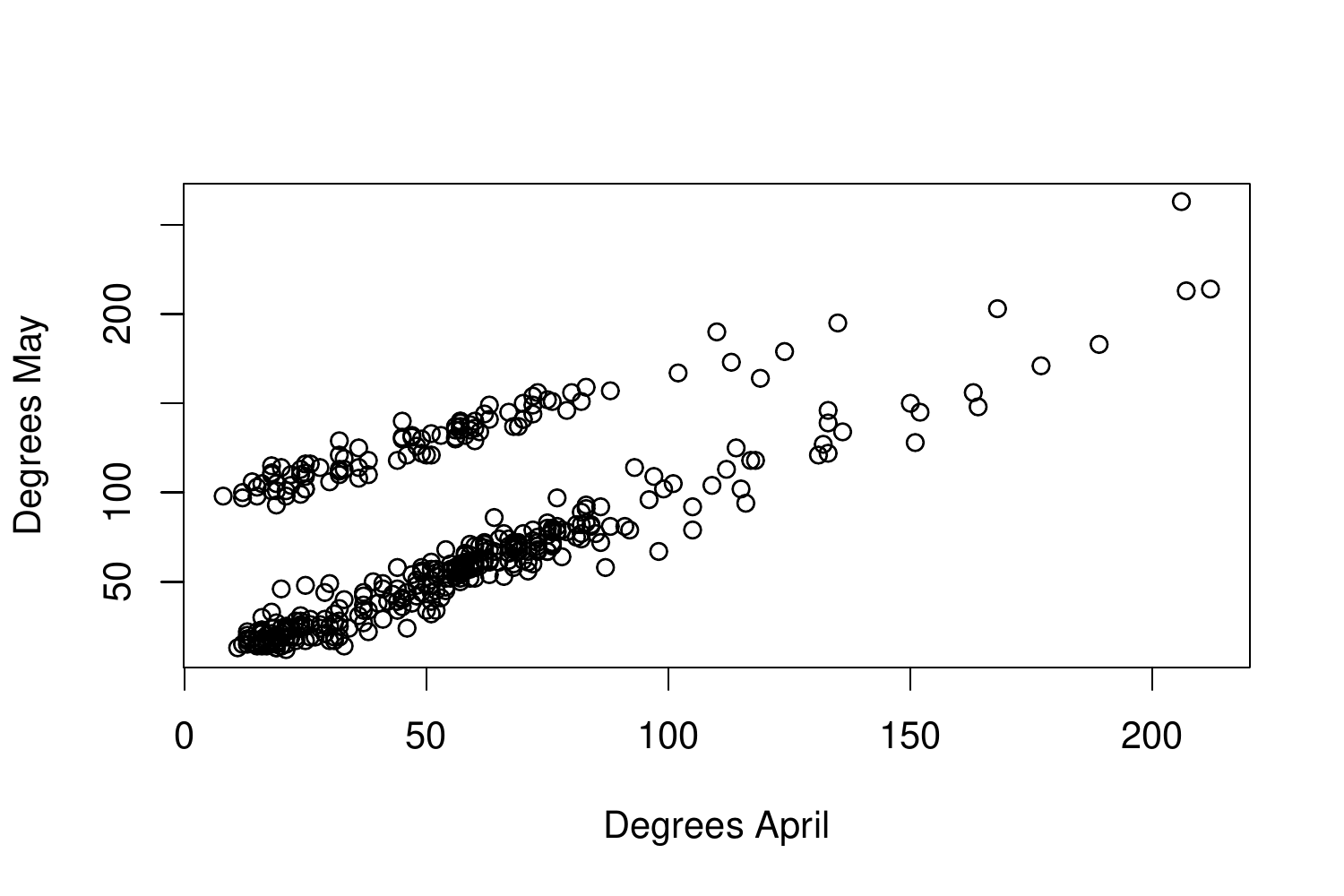}\\ 
	Vertices matched correctly & Degree Comparison, $\alpha=0.9$ 
\end{tabular}
	\caption{In the left panel, we plot the average matching accuracy ($\pm 1$ s.d.) of graph matching using \texttt{FAQ} initialized at $I_{431}$ when first choosing $100$ random vertices, denoted $\mathcal{V}$ from graph $B$ and then substituting $B[\mathcal V]\mapsfrom B[\mathcal V]+C$ (where $B$ is again binarized after noise is added) before matching and centering; here $C\sim \mathrm{ER}(100,q)$.
	Accuracy is plotted versus $q$
	In the right two panel, we plot the degrees of each vertex in April versus the degrees of the same vertex in May (with the $B[\mathcal V]\mapsfrom B[\mathcal V]+C$ substitution).}
	\label{fig:synthtwit}
\end{figure}

To explore the performance of the two approaches (centering and not centering) in the setting of different network topologies, we consider the following synthetic data experiment.
We choose $100$ random users from the twitter networks and add $C$ to their induced subgraphs in the May network, where $C\sim$ER$(100,q)$ \tb{(followed by again binarizing $B$); i.e., if the set of randomly chosen vertices is $\mathcal{V}$, then $B[\mathcal V]\mapsfrom B[\mathcal V]+C$ and the subsequent $B$ is binarized before matching or centering.}
We consider $q\in\{0.1,0.3,0.5,0.7,0.9\}$.
This experiment is simulating the setting where a fraction of the network changes its behavior from April to May; in this case by increasing their volume of mentions month to month.
For each value of $q$, we repeat this experiment $25$ times and plot the mean accuracy ($\pm$1 s.d.) of graph matching using \texttt{FAQ} initialized at $I_{431}$ both with and without USVT centering; see Figure \ref{fig:synthtwit}.
This experiment demonstrates the capacity for USVT to maintain $\hd$-matchability in the face of additive deviations in the network structure.  
These deviations have the effect of altering the graph topology month--to--month, and with enough signal, they have a precipitously negative impact on the performance of matching sans centering. 
Centering ameliorates this effect, and emphasizes the common signal in the networks by removing the effect of this additive noise.
\tb{It is interesting to note that for small values of $q$, centering negatively impacts algorithmic performance. 
We view this as potentially an artifact of the noise in these settings not being sufficient to obfuscate the true matching without centering.}

\begin{figure}[t!]
	\centering
	\includegraphics[width=.6\textwidth]{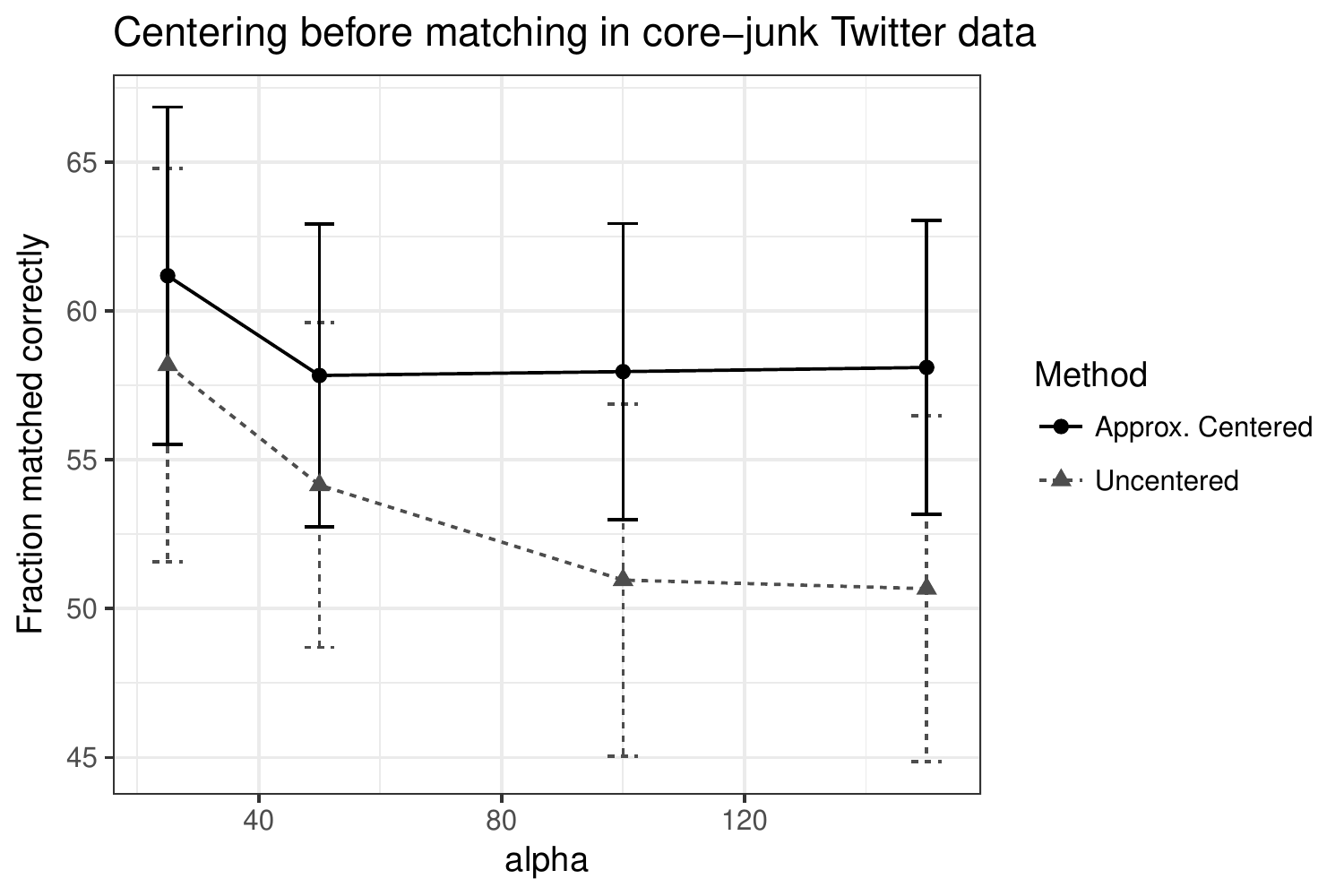} 
	\caption{We plot the average core matching accuracy ($\pm1$ s.d.) for $A[\mathcal{C}\cup\mathcal{J}_1]$ and $B[\mathcal{C}\cup\mathcal{J}_2]$ using \texttt{FAQ} initialized at $I_{n_c+n_j}$ (with $n_c=100$) against $n_j\in\{25,50,100,150\}$.
	Results are averaged over $100$ Monte Carlo iterates.}
	\label{fig:cjtwit}
\end{figure}

In the core-junk setting, the heterogeneity amongst the junk vertices offers a further setting for demonstrating the utility of USVT-centered graph matching.
To see this, we consider the following experiment:  choose $n_c=100$ uniformly sampled core vertices from the twitter network and $n_j\in\{25,50,100,150\}$ uniformly sampled junk vertices, $\mathcal{J}_1$, for the April graph and $n_j$ uniformly sampled junk vertices, $\mathcal{J}_2$, for the May graph.
As before, we match $A[\mathcal{C}\cup\mathcal{J}_1]$ and $B[\mathcal{C}\cup\mathcal{J}_2]$ using \texttt{FAQ} initialized at $I_{n_c+n_j}$ both with and without USVT centering; results are summarized in Figure \ref{fig:cjtwit}.
As seen previously, the ability of USVT-centering to ameliorate the degree/distributional heterogeneity (here amongst the junk vertices) leads to superior core label recovery compared to the uncentered matching setting.

\subsection{Connectomes}
\label{sec:connectome}

\tb{For our next example, we consider the the test-retest diffusion MRI data from \cite{landman2011multi}.
The dataset consists of test-retest pairs (used to evaluate reproducibility of the magnetization prepared rapid acquisition gradient echo (MPRAGE) image protocol).
Each scan is converted into a weighted connectome by considering 70 brain regions of interest (labeled according to the Desikan brain atlas \cite{desikan2006automated}) as the vertices, with edge weights counting the number of neural fiber bundles connecting the regions.
As vertices correspond to canonical brain regions of interest, it is natural to consider the true correspondence across graphs as being given by the identity mapping.}

\begin{figure}[t!]
	\centering
\begin{tabular}{cc}
	\includegraphics[width=.3\textwidth]{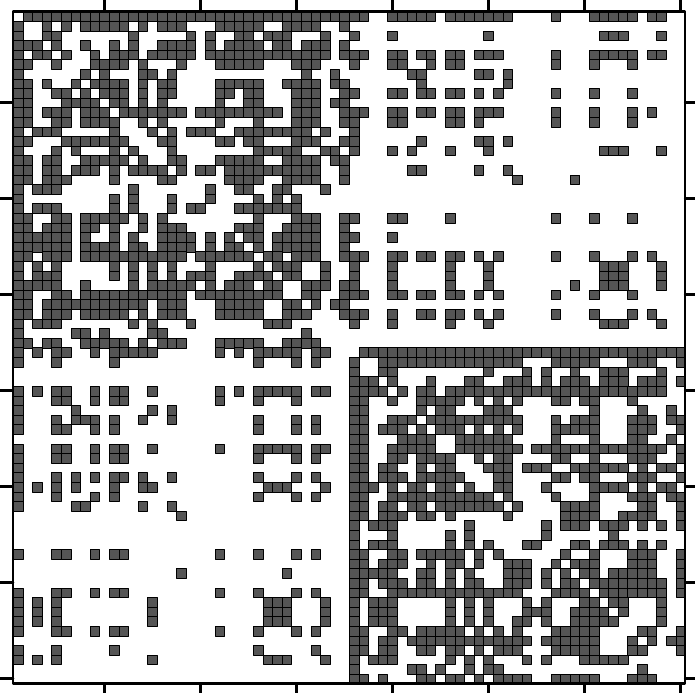}\hspace{10mm} &
	\includegraphics[width=.3\textwidth]{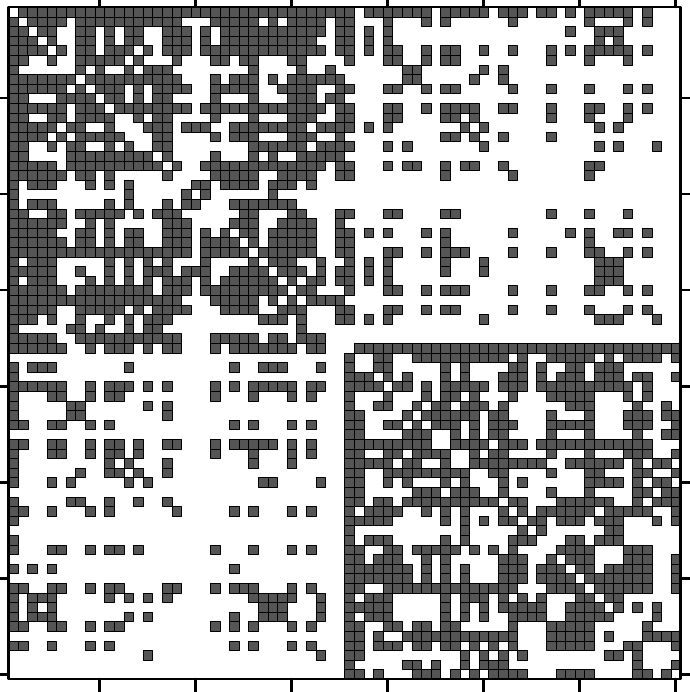}\\ 
	$G_1$& $G_2$ 
\end{tabular}
	\caption{Adjacency matrices of two sample brains from the dataset in \cite{landman2011multi}.}
	\label{fig:brains1}
\end{figure}

\tb{To illustrate the role of USVT centering in this data set, we first consider as an example a pair of graphs generated as above from the data in \cite{landman2011multi}.
The respective adjacency matrices for this graph pair are shown in Figure \ref{fig:brains1}.
Matching these brains directly using \texttt{FAQ} initialized at $I_{70}$ yields an estimated local optimum with 65 vertices correctly aligned across graphs; indeed, by permuting vertices $3,33,38,50,68$, we obtain a better objective function value than the GMP evaluate at $I_{70}$.
We seek to understand the ability of USVT centering, which is global in nature, to correct these local mismatches.}

\tb{To study this further, we apply a variant of the USVT procedure in which we automatically select the number of singular values to threshold by combining the ideas of USVT with the profile likelihood work of \cite{zhu2006automatic}; to wit, we select the threshold dimension via an elbow analysis of the SCREE plot of the singular values.
We chose this automated procedure rather than setting a singular value threshold because these graphs are weighted, and the common threshold of $2.01\sqrt{n}$ from \cite{chatterjee2015matrix,xu2017rates} is presented for the unweighted setting.
Centering the pair of graphs from Figure \ref{fig:brains1} recovers the identity $I_{70}$ as an estimated local minima of the GMP, and the global centering corrects the localized mismatch.}

\begin{figure}[t!]
	\centering
	\includegraphics[width=1\textwidth]{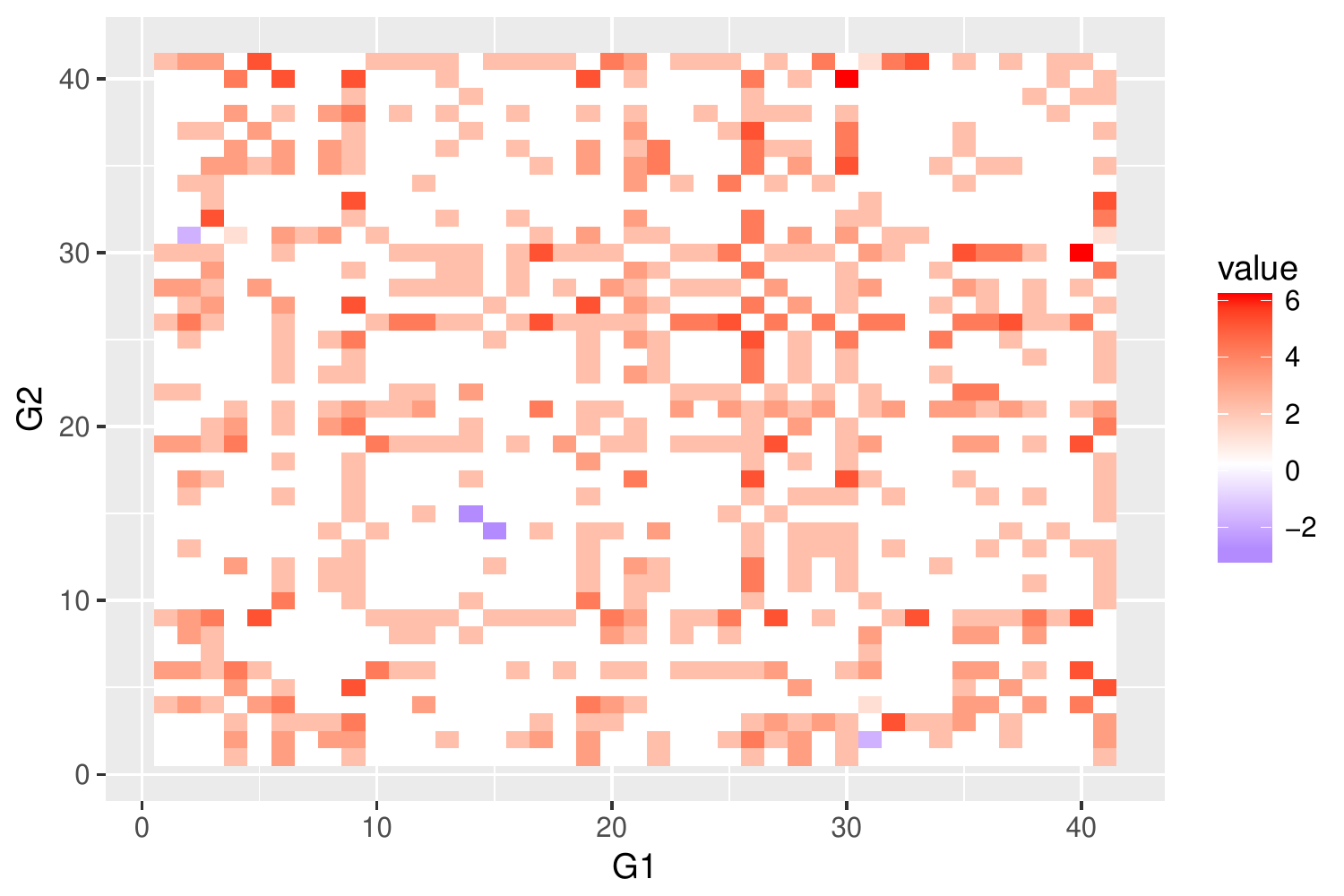}
	\caption{For the $\binom{41}{2}$ \emph{weighted} brain pairs, we plot a heatmap of the differences
$
D_c(i,j)-D(i,j)
$,
so that the $i,j$-th entry in the heatmap corresponds to the excess number of correct matches achieved by USVT centering.
Red values indicate more correctly matched via centering and blue values indicate more correctly matched via no centering. 
The color intensity indicates the value of $
D_c(i,j)-D(i,j)
$ achieved, with darker colors indicating more (in the red case) or less (in the blue case) vertices correctly matched after centering.  }
	\label{fig:brains2}
\end{figure}

\tb{Extending this to a 41 scan sample from \cite{landman2011multi} (hence, we consider 41 graphs each with 70 vertices), we run \texttt{FAQ} initialized at $I_{70}$ for the $\binom{41}{2}$ pairs of distinct graphs with both USVT centering and no centering.
When matching graph $i$ and graph $j$ for $\{i,j\}\in\binom{41}{2}$, we let
\begin{align*}
D_c(i,j)&=\#\text{ matched correctly by }\texttt{FAQ}\text{ initialized at }I_{70}\text{ in the USVT centered case}\\
D(i,j)&=\#\text{ matched correctly by }\texttt{FAQ}\text{ initialized at }I_{70}\text{ in the un centered case}.
\end{align*}
In Figure \ref{fig:brains2}, we plot a heatmap of the $\binom{41}{2}$ differences
$
D_c(i,j)-D(i,j)
$,
so that the $i,j$-th entry in the heatmap corresponds to the excess number of correct matches achieved by USVT centering.
Red values indicate more correctly matched via centering and blue values indicate more correctly matched via no centering. 
The color intensity indicates the value of $
D_c(i,j)-D(i,j)
$ achieved, with darker colors indicating more (in the red case) or less (in the blue case) vertices correctly matched after centering. 
The figure demonstrates that the phenomena observed in the graphs in Figure \ref{fig:brains1} was not an anomaly.
Only $2$ pairs see an improvement in matching accuracy when not centering, while $284$ pairs see an improvement in matching accuracy when USVT centering.
Moreover, while many of the mismatches are local in nature, they are nonetheless ameliorated by the global USVT centering procedure.}

\begin{figure}[t!]
	\centering
\begin{tabular}{cc}
	\includegraphics[width=.5\textwidth]{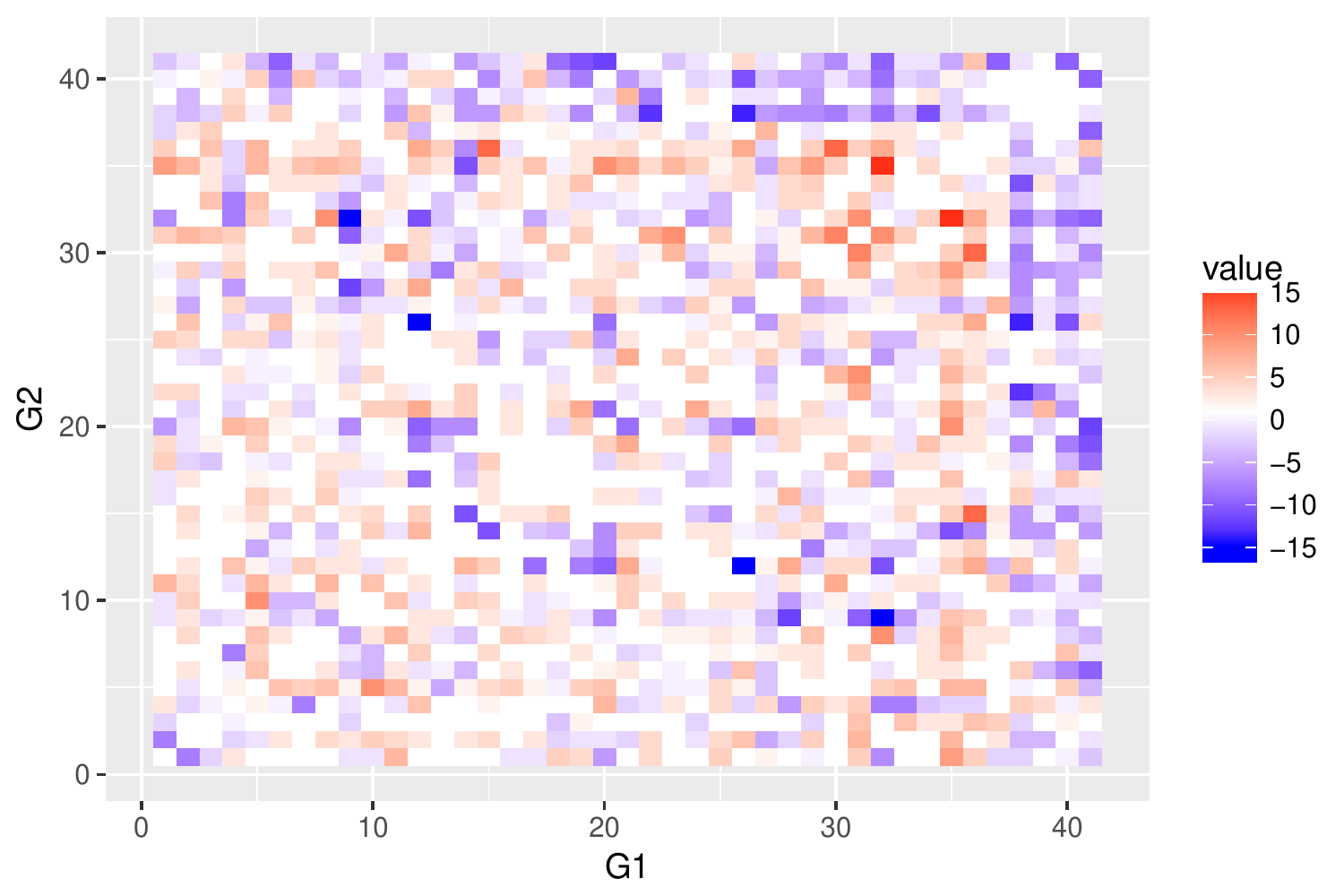} &
	\includegraphics[width=.5\textwidth]{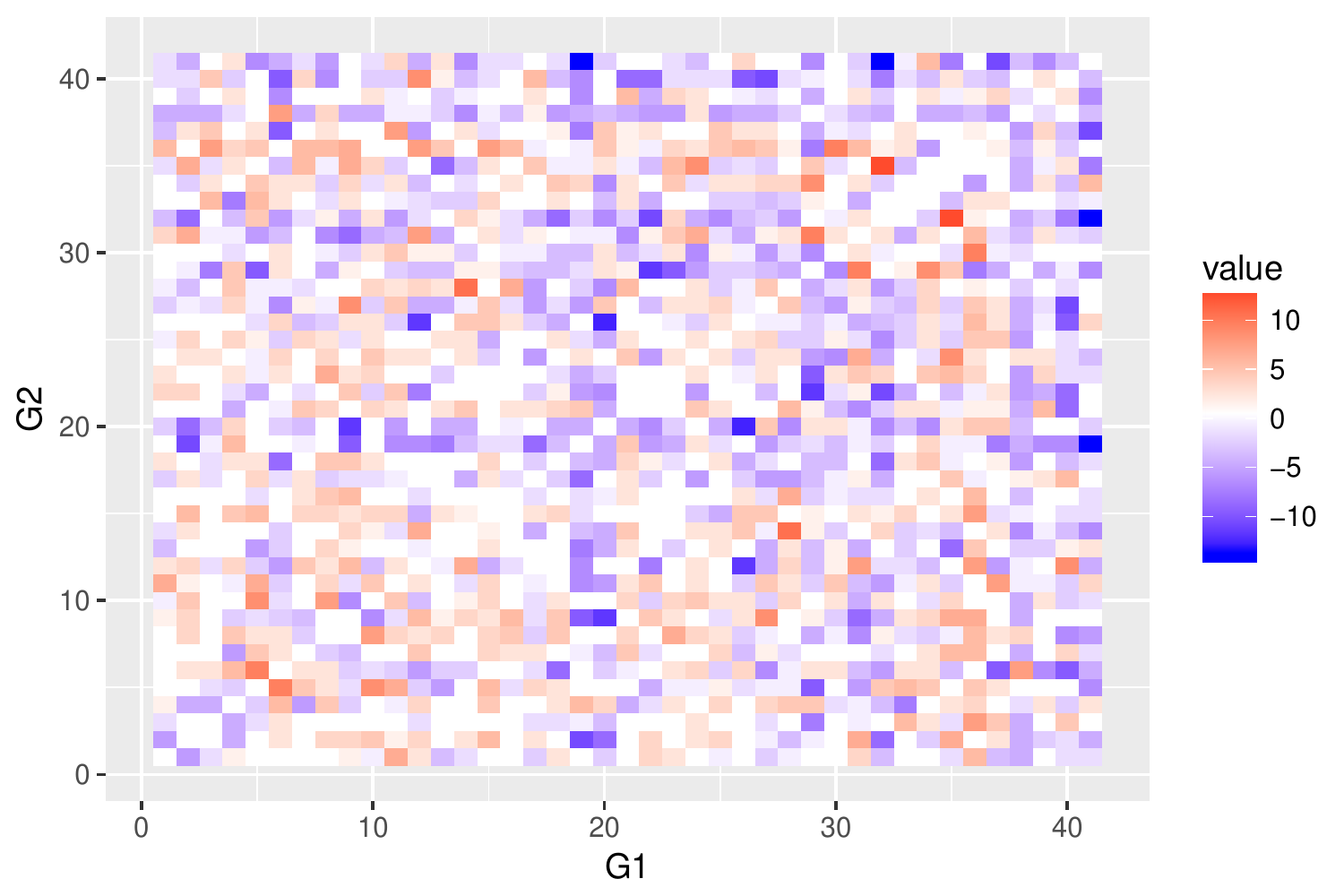}\\ 
	USVT with $t=2\sqrt{70}$& USVT with $t=\sqrt{70}$
\end{tabular}
	\caption{For the $\binom{41}{2}$ \emph{unweighted} brain pairs, we plot a heatmap of the differences
$
D_c(i,j)-D(i,j)
$,
so that the $i,j$-th entry in the heatmap corresponds to the excess number of correct matches achieved by USVT centering.
Red values indicate more correctly matched via centering and blue values indicate more correctly matched via no centering. 
The color intensity indicates the value of $
D_c(i,j)-D(i,j)
$ achieved, with darker colors indicating more (in the red case) or less (in the blue case) vertices correctly matched after centering.  
In the left panel, we center by USVT with $t=2\sqrt{70}$; in the right panel by $t=\sqrt{70}$.}
\end{figure}

\tb{If we consider running the same experiment on the unweighted brain graphs (using USVT centering with threshold $t$), we see the delicate nature of the USVT threshold in data applications.
We note here that, again, in this unweighted case, we do not zero out the diagonal of $\hQ_i$ in the USVT step in this real data example.
When $t=2\sqrt{70}$, $267$ pairs achieved improved matching performance when USVT centering first, 
and $271$ pairs achieved improved matching performance when not centering.
When $t=\sqrt{70}$, $264$ pairs achieved improved matching performance when USVT centering first, 
and $343$ pairs achieved improved matching performance when not centering.
These results suggest two important take-aways:
First, performance is intimately tied to properly thresholding; and second, in this example, USVT-centering is more effective in the weighted edge case.
This suggests both that the magnitude of weights are contributing significantly to the mismatch and that USVT centering is effective at ameliorating this edge weight heterogeneity; this is not entirely unexpected as the centering is precisely trying to eliminate the different edge probability/weight structures across the $Q_i$.}


\section{Discussion} 
\label{sec:discussion}
Understanding the limits of $\delta$-matchability is an essential step in robust multiple graph inference regimes.
When graphs are not $\delta$-matchable---i.e., the true node correspondence cannot be recovered in the face of noise---paired graph inference methodologies that utilize the across graph correspondence (see, for example, \cite{MT2,asta2014geometric}) cannot gainfully be employed.
Non-$\delta$-matchability can limit analysis to methods which rely on graph statistics which are invariant to relabeling of the vertices, which can be useful but lack the full power of their parametric (with the labeling as a parameter) counterparts.
In this paper, we establish initial theoretical results on $\delta$-matchability when the graphs to be match differ in distribution and when only a fraction of the graphs are matchable.

While our theoretical results and subsequent simulations and experiments provide a basis for a deeper understanding of the effect that distributional heterogeneity has on $\delta$-matchability, there is still much to be done.
For example, in our present USVT centering step, the across graph correlation provided by $R$ is not utilized. 
\tb{We suspect that the error in the USVT steps could be greatly reduced by leveraging $R$ or an estimate thereof.}
We are also exploring graph normalization strategies other than centering, such as kernel smoothing (using a small number of a priori known correspondences to choose the proper smoothing kernels), which may be more appropriate in the presence of multiplicative or other non-linear noise structures.
We suspect that the growth rate of $k_0$ in Theorem \ref{thm:estimatedcenteredmatched} is not sharp, as in simulation and real data settings the true correspondence is almost perfectly recovered via USVT centering before matching; however, sharpening the lower bound on $k_0$ with our present methods does not seem feasible and new ideas and techniques need be employed.

We suspect that the bounds on $n_j$ versus $n_c$ obtained in Theorem \ref{thm:coreestcenteredmatched}
are also suboptimal.
As an informal argument we can consider known results for the quadratic assignment problem for i.i.d. entries \cite{cela2013quadratic}.
In  the uncorrelated dense homogeneous Erd\H os-R\'enyi setting, these results imply that the best best solution to the graph matching problem, Eq.~\eqref{eq:GMP} will reduce the objective function $\Theta(n_j^{3/2})$ as compared to a random guess.
On the other hand, in the dense homogeneous Erd\H os-R\'enyi setting with constant correlation, the best solution is $\Theta(n_c^2)$ better than a random guess.
Under the heuristic that in order for the core vertices to be matched correctly, the signal from correlation must be greater than the possible improvements in the all noise setting, we can conjecture that a core matchability threshold at approximately $n_c=O(n^{3/4})$ may be possible. 
Using our present proof technique, we are unable to achieve this rate.
This heuristic argument, though problematic, provides a potential guidepost for future work.

In addition, while the theoretical results presented herein are for the case that the graphs are simple undirected graph with no edge-weights, the graph matching framework of Definition \ref{def:GMP} is flexible, allowing us to accommodate many of the features---both those considered above, and additional eccentricities---inherent to real data settings.
In the weighted, loopy setting, the matching can occur between the weighted adjacency matrices or normalized Laplacian matrices, $\mathcal L_A=D_A^{-1/2}AD_A^{-1/2}$ (where $D_A$ is the diagonal matrix with $(i,i)$-th entry equal to $\sum_j D_A(i,j)$) and the similarly defined $\mathcal{L}_B$.
To match directed graphs, the graphs can either be made undirected (for example, by matching $A'=(A+A^T)/2$ to $B'=(B+B^T)/2$) or the directed adjacency matrices can be directly plugged into Eq.~\eqref{eq:GMP}.
Developing similar results to Theorems \ref{thm:estimatedcenteredmatched} and \ref{thm:coreestcenteredmatched} in models (akin to our CorrER model) that incorporate these graph features as well as additional vertex and edge features is a natural next step.

The information and computational limits for $\delta$-matchability are still open problems for which we have pushed the boundaries, but significant more work is to be done. 
These problems are in analogy to the recently addressed problems of detection and recovery for the planted partition and planted clique problems for a single graph \cite{Feldman2013-we,Mossel2014-lk,Abbe2015-ua}.
For these settings, exact fundamental limits have been established and polynomial time algorithms have been shown to achieve or nearly achieve these limits.
Obtaining similar results for the graph matching problem are key steps towards a robust statistical framework for multiple graph inference.


\subsection*{Acknowledgments}

Dr. Sussman's contribution to this work was partially supported by a grant from MIT Lincoln Labs and the Department of Defense.
This material is based on research sponsored by the Air Force Research Laboratory and DARPA under agreement number FA8750-18-2-0066. The U.S.
Government is authorized to reproduce and distribute reprints for Governmental purposes notwithstanding any copyright notation thereon. The views
and conclusions contained herein are those of the authors and should not
be interpreted as necessarily representing the official policies or endorsements, either expressed or implied, of the Air Force Research Laboratory
and DARPA or the U.S. Government.
We would also like to thank Prof. Carey E. Priebe, Prof. Minh Tang, and Joshua Cape for their helpful discussions in the writing of this manuscript.


\appendix

\section{Proofs}
Herein we collect proofs of the main theoretical results in this manuscript.  Before stating our proofs, we first state some well-known facts about the bivariate Bernoulli distribution.
\label{sub:bibern}
\tb{Indeed, if $(A,B)\sim\HER$ model, then for each $\{u,v\}\in\binom{[n]}{2}$, $A(u,v)$ and $B(u,v)$ can be realized as a bivariate Bernoulli random variable.}
This will be a key insight in the proof of our main results, Theorems \ref{thm:centeredmatched} and \ref{thm:estimatedcenteredmatched}.

Spelling this out further, a pair of Bernoulli random variables $(X,Y)$ has a BiBernoulli distribution with 
\begin{equation}
     (X,Y)\sim \mathrm{BiBern} \begin{pmatrix}
    p_{11} & p_{10} \\ p_{01} & p_{00}
    \end{pmatrix}
\end{equation} 
if $\Pr[X=x,Y=y] = p_{xy}$ for each $(x,y)\in \{0,1\}^2$.
A key property of BiBernoulli random variables is that they can be generated by a triple of independent Bernoulli random variables.
For $(X,Y)$ as above, setting $Z_0\sim \mathrm{Bern}(p_{11}+p_{10}$), $Z_1\sim \mathrm{Bern}(p_{01}/(p_{01}+p_{00}))$, and $Z_2\sim \mathrm{Bern}(p_{11}/(p_{11}+p_{10}))$, with $Z_0,Z_1$ and $Z_2$ independent yields
    $$(X,Y)\stackrel{\mathcal{L}}{=}(Z_0,Z_0Z_2+(1-Z_0)Z_1).$$ 
\tb{In the $(A,B)\sim\HER$ model, we then have that
 $$(A(u,v),B(u,v))\stackrel{\mathcal{L}}{=}(Z_0(u,v),Z_0(u,v)Z_2(u,v)+(1-Z_0(u,v))Z_1(u,v)),$$
 where
 \begin{align*}
Z_0(u,v)&\sim \mathrm{Bern}\left(Q_1(u,v)
\right)\\
Z_1(u,v)&\sim \mathrm{Bern}\left(
\frac{Q_2(u,v)(1-Q_1(u,v))-\text{Cov}(A(u,v),B(u,v))
}{1-Q_1(u,v)}
\right)\\
Z_2(u,v)&\sim \mathrm{Bern}
\left(
\frac{Q_2(u,v)Q_1(u,v)+\text{Cov}(A(u,v),B(u,v))
}{Q_1(u,v)}
\right)
 \end{align*}
 are independent Bernoulli random variables.}


\subsection{Proof of Theorem \ref{thm:centeredmatched}}
\label{AP:TH5}
The key to the Proof of Theorem \ref{thm:centeredmatched} is the well-known McDiarmid's inequality \cite{mcdiarmid1989method}.
\tb{\begin{proposition}[McDiarmid's inequality]
\label{prop:kim} 
Let $X_1,\dotsc,X_n\in\mathcal{X}$ be a sequence of independent random variables. 
Let $f:\mathcal{X}^n\mapsto \Re$ be such that 
for all $i\in[n]$, and for all $x_1,\ldots,x_n,x_i'\in\mathcal{X}$ 
$$ 
   |f(x_1,\dotsc, x_{i-1},x_i,x_{i+1},\dotsc,x_n) - f(x_1,\dotsc,x_{i-1},x_i',x_{i+1},\dotsc,x_n)|\leq c_i.
$$
If $Y=f(X_1,\dotsc,X_n)$, then for any $t>0$,
$$\Pr[|Y-\Ex[Y]| \geq t ] \leq 2 \text{exp}\left\{-\frac{-2t^2}
{\sum_{i=1}^n c_i^2}
\right\}.$$
\end{proposition}}

\noindent{\it Proof of Theorem \ref{thm:centeredmatched}:} 
Let $G_1,G_2\sim\HER$, and consider $\tA=A-\e(B)$ and $\tB=B-\e(B)$.

 If $X$ and $Y$ are $\varrho$-correlated Bernoulli random variables with respective parameters $p$ and $q$, then
it follows that
\begin{align*}
\e[(X-p)(Y-q)]=\text{Cov}(X,Y)=\varrho\sqrt{pq(1-p)(1-q)}.
\end{align*}
Let $P\in\Pi(n)$ be a permutation matrix that permutes exactly $k$ labels and let $\tau$ be the associated permutation for $P$. 
Note that 
\tb{$$\mathcal{E}_P=\left\{ \{u,v\}\in\binom{V}{2}\text{ s.t. } \{\tau(u),\tau(v)\}\neq\{u,v\}\right\}$$
satisfies
$$|\mathcal{E}_P|=\underbrace{\binom{k}{2}-T_P}_{
\begin{tabular}{@{}l@{}}pairs in $\mathcal{E}_P$ where both \\u and v are perm. by $\tau$\end{tabular}}+\underbrace{(n-k)k}_{
\begin{tabular}{@{}l@{}}
pairs in $\mathcal{E}_P$ where only one \\of u or v are perm. by $\tau$\end{tabular}},$$
where 
$$
T_P=\big\{\{u,v\}\text{ s.t. }u=\tau(v),\,v=\tau(u)\big\}
$$ 
is the number of transpositions induced by $P$.
Note that $|T_P|\leq k/2$,
and so
\begin{align}
\label{eq:epnd}
|\epsilon_P|\geq nk-\frac{k^2}{2}-\frac{k}{2}=k(n-1-\frac{k}{2}).
\end{align}
For each $Q\in\Pi(n)$, define the matrix $\ck_{A,B,Q}$ as in Eq. \ref{eq:K}. 
We have then
\begin{align}
\notag
\frac{1}{2}\,\,\e\,\left(\text{tr}\left(\tA\tB\right)-\text{tr}\left(\tA P\tB P^T\right)\right)
&=\sum_{\{u,v\}\in\binom{V}{2}}\text{Cov}(A(u,v),B(u,v))-\text{Cov}(A(u,v),B(\tau(u),\tau(v)))\\
\notag
&=\sum_{\{u,v\}\in\mathcal{E}_P}\text{Cov}(A(u,v),B(u,v))-\text{Cov}(A(u,v),B(\tau(u),\tau(v)))\\
\label{eq:eps1}
&=\sum_{\{u,v\}\in\mathcal{E}_P}\text{Cov}(A(u,v),B(u,v))\\
\notag
&=\frac{1}{2}\left[\text{tr}\left(\ck_{A,B,I_n}\right)-\text{tr}\left(\ck_{A,B,P}\right)\right].
\end{align}
For ease of notation, we define 
\begin{align*}
X_P&:=\frac{1}{2}\left(\text{tr}(\tA\tB)-\text{tr}(\tA P\tB P^T)\right),
\end{align*}
so that
\begin{align*}
\mathcal{X}_P&=\e(X_P)=\frac{1}{2}\left[\text{tr}\left(\ck_{A,B,I_n}\right)-\text{tr}\left(\ck_{A,B,P}\right)\right].
\end{align*}
As $X_P=\sum_{\{u,v\}\in\mathcal{E}_P} \tA(u,v)\left( \tB(u,v)-\tB(\tau(u),\tau(v))\right)$, and each $\{A(u,v),B(u,v)\}$ pair is a function of three independent Bernoulli random variables, we have that
$X_P$ is a function of 
$$\left\{\{A(u,v),B(u,v)\}\text{ s.t. }\{u,v\}\in\mathcal{E}_P \right\}$$
and so is a function of $N_P:=3|\mathcal{E}_P|$ independent Bernoulli random variables,
where $N_P$ then satisfies
\begin{align*}
N_P\leq 3\left(nk-\frac{k^2}{2}\right).
\end{align*}}

\tb{Changing the value of one of these Bernoulli random variables leaving all others fixed can change the value of at most one $\{A(u,v),B(u,v)\}$ pair, and this pair appears in two terms in the sum $X_P$.
As each term in the sum of $X_P$ is bounded in $[-2,2]$, we have that, in the notation of Proposition \ref{prop:kim}, each $c_i$ can be uniformly set to $8$.
Proposition \ref{prop:kim} then yields (setting $t=\mathcal{X}_P$)
\begin{align}
\label{eq:lowerbnd}
\Pr\left(X_P\leq 0\right)&\leq\Pr\left(|X_P-\mathcal{X}_P|\geq \mathcal{X}_P\right)\\
&\leq 2\exp\left\{-C\,\,\frac{\mathcal{X}_P^2}{N_P}\right\}\notag\\
&\leq 2\exp\left\{-C\,\, \frac{\mathcal{X}_P^2}{nk-\frac{k^2}{2}}\right\}\notag
\end{align}
with $C$ being an appropriate positive constant that may change line to line.
If $\mathcal{X}_P=\omega(k\sqrt{n\log(n)})$ then
\begin{align}
\label{eq:bndbndbnd}
\Pr\left(X_P\leq 0\right)&\leq 2\exp\left\{-\omega(k\log n)\right\}
\end{align}}

To finish the proof, we apply a union bound on all such $P$.
The number of such permutations $P$ that permute $k$ vertex labels is upper bounded by \tb{$n^{k}$.}  Combining with Eq. (\ref{eq:bndbndbnd}), we have that under the assumptions of the Theorem
\begin{align*}
\Pr[ \min_{P\in \Pi(n)\setminus \{I_n\}}& X_P \leq 0 ]\leq \sum_{k=2}^{n} 2n^{k}\exp\left\{-\omega(k\log n)\right\}\\
=&\sum_{k=2}^{n} 2\exp\left\{-\omega(k\log n)\right\}\\
=& 2\exp\left\{-\omega(\log n)\right\},
\end{align*}
as desired.


\subsection{Proof of Theorem \ref{thm:estimatedcenteredmatched}}
\label{pf:estimatedcenteredmatched}

\tb{Key to the proof of Theorem \ref{thm:estimatedcenteredmatched} are the following Lemmas, adapted here from \cite[Lemmas 1 and 2]{xu2017rates}.
\begin{lemma}
\label{lem:lem1}
Let $X,Y\in\mathbb{R}^{n\times n}$. 
Suppose $\tau\geq (1+\delta)\|X-Y\|$ for $\delta>0$ fixed.
Let $X=\sum_i \sigma_i(X)u_i v_i^T$ be the singular value decomposition of $X$, and let 
$$\breve X=\sum_{i:\sigma_i(X)>\tau} \sigma_i(X)u_i v_i^T.$$
Then
$$\|\breve X-Y\|_F^2\leq 16\min_{0\leq r\leq n}\left(\tau^2 r+\left(\frac{1+\delta}{\delta} \right)^2\sum_{i\geq r+1}\sigma_i^2(Y)
\right),$$
where $\sigma_i(Y)\geq \sigma_i(Y)\geq\cdots\geq \sigma_i(Y)$ are the singular values of $Y$.
\end{lemma}}

\tb{\begin{lemma}
\label{lem:lem2}
Let $X\sim$ER($Q$); i.e., $X$ is a hollow, symmetric matrix with $X(u,v)\stackrel{ind.}{\sim}Q(u,v)$.
Assume $Q\leq r$ for some $r>0$.
If $nr\geq C\log n$ for a constant $C>0$, then for all $c>0$ there exists a constant $\kappa>0$ such that 
$$\p(\|X-\e(X)\|\leq \kappa\sqrt{nr})\geq 1-n^{-c}.$$
\end{lemma}}

\tb{Note first that 
\begin{align*}
\| A-Q_1\|&\leq \| A-\e(A)\|+\| Q_1-\e(A)\|\\
&\leq \| A-\e(A)\|+r_1,
\end{align*}
with the analogous result holding for $B$.
Under the assumptions of Theorem \ref{thm:estimatedcenteredmatched}, Lemma \ref{lem:lem2} implies that there exists constants $c_1,c_2>0$ such that 
$$\| A-\e(A)\|\leq c_1\sqrt{nr_1}\text{ and }\| B-\e(B)\|\leq c_2\sqrt{nr_2} $$
with probability at least $1-n^{-2}$, and therefore there exists constants 
$c'_1,c'_2>0$ such that 
$$\| A-Q_1\|\leq c'_1\sqrt{nr_1}\text{ and }\| B-Q_2\|\leq c'_2\sqrt{nr_2} $$
with probability at least $1-n^{-2}$.}

\tb{Next, we apply Lemma \ref{lem:lem1} with $X=A$, $Y=Q_1$ (resp., $X=B$, $Y=Q_2$).
With probability at least $1-n^{-2}$,
there exists constants $a_1,C_1>0$ such that if $\tau=a_1\sqrt{nr_1}$ then (where $\breve A$ is as defined in the USVT pseudocode, Algorithm \ref{alg:USVT})
\begin{align*}
\|\breve A-Q_1\|_F^2\leq 16\left(a_1^2 nr_1 d_1+C_1\sum_{i\geq d_1+1}\sigma_i^2(Q_1)
\right)=O(nr_1 d_1),
\end{align*}
where the equality follows from the rank assumption on $Q_1$ in Theorem \ref{thm:estimatedcenteredmatched}; similarly, for $B$ we have
\begin{align*}
\|\breve B-Q_2\|_F^2=O(nr_2 d_2).
\end{align*}}

\tb{Combining the above, we have that there exists an event $\mathcal{E}_1$ such that $\p(\mathcal{E}_1)\geq 1-n^{-2}$, and on $\mathcal{E}_1$,
\begin{align*}
\|\widehat Q_1-\e(A)\|_F^2&\leq \|\breve A-Q_1\|_F^2=O(nr_1 d_1);\\
\|\widehat Q_2-\e(B)\|_F^2&\leq \|\breve B-Q_2\|_F^2=O(nr_2 d_2).
\end{align*}}

\tb{To prove Theorem \ref{thm:estimatedcenteredmatched}, we proceed as follows.
Fix $P\in\Pi(n)$ so that $\sum_i P(i,i)=n-k$; i.e., $P$ permutes exactly $k$ labels.
Two simple applications of the triangle inequality yields that
\begin{align*}
\|\hA-P\hB P^T\|_F
&=\|A-\hQ_1-PB P^T+P\hQ_2P^T\|_F
\\
&=\|\tA-(\hQ_1-\e(A))-P\tB P^T-P(\e(B)-\hQ_2)P^T\|_F\\
&\geq \|\tA-P\tB P^T\|_F-\|\hQ_1-\e(A)\|_F-\|\e(B)-\hQ_2\|_F,
\end{align*}
and 
\begin{align*}
\|\hA-\hB \|_F
&\leq \|\tA-\tB \|_F+\|\hQ_1-\e(A)\|_F+\|\e(B)-\hQ_2\|_F.
\end{align*}
Combining the above, we have that
\begin{align*}
\|\hA-P\hB P^T\|_F-\|\hA-\hB \|_F
&\geq \|\tA-P\tB P^T\|_F-\|\tA-\tB \|_F-2\|\hQ_1-\e(A)\|_F-2\|\e(B)-\hQ_2\|_F.
\end{align*}
In the proof of Theorem \ref{thm:centeredmatched}, if we set $t=\mathcal{X}_P/2$ in Eq. \ref{eq:lowerbnd} when applying McDiarmid's inequality, then under the assumptions of Theorem \ref{thm:estimatedcenteredmatched}, there exists an event $\mathcal{E}_{2,P}$ with $\p(\mathcal{E}_{2,P})\geq 1-2\text{exp}\{-\omega(k\log n)\}$ such that on $\mathcal{E}_{2,P}$, 
$$ \|\tA-P\tB P^T\|^2_F-\|\tA-\tB \|^2_F =\omega\left(kn^{1/2+\delta}\sqrt{\log n}\right).
$$}

\tb{Next, note that 
\begin{align*}
\e(\|\tA\|^2_F)&= 2\sum_{\ell\geq j}Q_1(\ell,j)(1-Q_1(\ell,j))<2r_1\binom{n}{2}, \text{ and}\\\
\e(\|\tB \|^2_F)&=2\sum_{\ell\geq j}Q_2(\ell,j)(1-Q_2(\ell,j))<2r_2\binom{n}{2}.
\end{align*}
Hoeffding's inequality (see, for example, \cite{fcci}) yields that 
\begin{align}
\label{eq:h1}
\|\tA\|^2_F<3r_1\binom{n}{2}\Rightarrow \|\tA\|_F<2\sqrt{r_1}n;\\
\label{eq:h2}
\|\tB\|^2_F<3r_2\binom{n}{2}\Rightarrow \|\tB\|_F<2\sqrt{r_2}n.
\end{align}
with probability at least 
\begin{align}
\label{eq:phoef}
\p(\mathcal{E}_3)\geq1-&\text{exp}\left\{-\frac{2r_1^2\binom{n}{2}}{(1-2r_1)^2}\right\}-\text{exp}\left\{-\frac{2r_2^2\binom{n}{2}}{(1-2r_2)^2}\right\}\geq 1-2e^{-\omega(\log^2 n)}
\end{align}
(where the last inequality followed from the assumptions in the Theorem, as under the assumptions $nr_i=\omega(\log n)$).
Therefore there exists an event $\mathcal{E}_3$
such that Eq. (\ref{eq:h1}) and (\ref{eq:h2}) hold on $\mathcal{E}_3$, and $\p(\mathcal{E}_3)\geq 1-2e^{-\omega(\log^2 n)}.$}

\tb{Writing 
\begin{align*}
\|\tA-P\tB P^T\|^2_F-\|\tA-\tB \|^2_F &=(\|\tA-P\tB P^T\|_F-\|\tA-\tB \|_F)(\|\tA-P\tB P^T\|_F+\|\tA-\tB \|_F)\\
&\leq(\|\tA-P\tB P^T\|_F-\|\tA-\tB \|_F)(2\|\tA\|_F+2\|\tB\|_F),
\end{align*}
we see that
\begin{align*}
\|\hA-P\hB P^T\|_F-\|\hA-\hB \|_F\geq \frac{\|\tA-P\tB P^T\|^2_F-\|\tA-\tB \|^2_F}{2\|\tA\|_F+2\|\tB\|_F}-2\|\hQ_1-\e(A)\|_F-2\|\e(B)-\hQ_2\|_F.
\end{align*}
We see then that on $\mathcal{E}_1\cap\mathcal{E}_{2,P}\cap\mathcal{E}_3$,
\begin{align}
\label{eq:k0}
\|\hA-P\hB P^T\|_F-\|\hA-\hB \|_F\geq \frac{\omega\left(kn^{1/2+\delta}\sqrt{\log n}\right)}{4n(\sqrt{r_1}+\sqrt{r_2})}-O(\sqrt{nr_1d_1})-O(\sqrt{nr_2d_2}).
\end{align}
Let $k_0$ be such that 
$$k_0=\omega\left(
\frac{n^{1-\delta}}{\sqrt{\log n}}(\sqrt{r_1}+\sqrt{r_2})\left(\sqrt{r_1 d_1}+\sqrt{r_2 d_2}\right) \right).$$
If $P\in\Pi(n,k)$ for $k\geq k_0$, 
$\|\hA-P\hB P^T\|_F-\|\hA-\hB \|_F>0$ on
$\mathcal{E}_1\cap\mathcal{E}_{2,P}\cap\mathcal{E}_3$
where 
\begin{align*}
\p(\mathcal{E}_1\cap\mathcal{E}_{2,P}\cap\mathcal{E}_3)&=1-\p(\mathcal{E}_1^c\cup\mathcal{E}_{2,P}^c\cup\mathcal{E}_3^c)\\
&\geq 1-n^{-2}-2e^{-\omega(k\log n)}-2e^{-\omega(\log^2n)}.
\end{align*}}

\tb{Define the event 
$$\widehat{\mathcal{E}}_P=\left\{ \|\hA-P\hB P^T\|_F\leq\|\hA-\hB \|_F \right\}$$
and note that
$$\widehat{\mathcal{E}}_P\subset \mathcal{E}_1^c\cup\mathcal{E}_{2,P}^c\cup\mathcal{E}_3^c.$$
Combined, this yields
\begin{align*}
\p&\left(\bigcup_{k\geq k_0}\bigcup_{P\in\Pi(n,k)}\widehat{\mathcal{E}}_P\right)=\p\left(\exists P\in\Pi(n,k)\text{ with }k\geq k_0\text{ and }\|\hA-P\hB P^T\|_F\leq\|\hA-\hB \|_F\right)\\
&\leq \p\left(\bigcup_{k\geq k_0}\bigcup_{P\in\Pi(n,k)}\mathcal{E}_1^c\cup\mathcal{E}_{2,P}^c\cup\mathcal{E}_3^c\right)=\p\left(\mathcal{E}_1^c\cup\left(\bigcup_{k\geq k_0}\bigcup_{P\in\Pi(n,k)}\mathcal{E}_{2,P}^c\right)\cup\mathcal{E}_3^c\right)\\
&\leq \p(\mathcal{E}_1^c)+\sum_{k\geq k_0}\sum_{P\in\Pi(n,k)}\p(\mathcal{E}_{2,P}^c)+\p(\mathcal{E}_3^c)\\
&\leq n^{-2}+\sum_{k\geq k_0}\sum_{P\in\Pi(n,k)}2e^{-\omega(k\log n)}+2e^{-\omega(\log^2n)}\\
&\leq n^{-2}+\sum_{k\geq k_0}n^k2e^{-\omega(k\log n)}+2e^{-\log^2n}\\
&\leq 5n^{-2},
\end{align*}
as desired.  }


\subsection{Proof of Theorem \ref{thm:matchedcorescentered}}
\label{AP:stochmatchedcores}
For a given permutation $\tau$ on $[n]$, we define the permutation $\tau_{\id}$ uniquely as follows:
\begin{equation}
  \tau_{\id}(i)= \begin{cases}
    i, &\text{ if }i\in \mathcal{C} \\
    \tau^k(i), &\text{ if }i\in \mathcal{J}, \text{ where } k=\min\{\ell\geq 1: \tau^\ell(i)\in \mathcal{J}\}.
  \end{cases}\label{eq:sigma_id}
\end{equation}
For example, if $n_c=4$, $n_j=4$, and 
$$\tau=\begin{pmatrix}1&2&3&4&5&6&7&8\\
3&6&7&1&8&5&4&2  \end{pmatrix},$$
then
$$\tau_{\id}=\begin{pmatrix}1&2&3&4&5&6&7&8\\
1&2&3&4&8&5&7&6  \end{pmatrix}.$$
For a permutation matrix $P$, we define $P_{\id}\in \mathcal{P}_n^\mathcal{C}$ analogously, where we recall here that $\pnc$ is the set of permutation matrices $P$ in $\Pi(n)$ satisfying $\sum_{i=1}^{n_c}P(i,i)=n_c$ (i.e., fixing all core labels).
Define $\mathcal{O}_n^\mathcal{C}=\Pi_n\setminus \pnc$.
Define the events
\begin{align*}
\mathcal{B}_1&:=\{\exists Q\in \mathcal{O}_n^\mathcal{C}\text{ s.t. } \forall P\in \pnc,\,\|\tA Q-Q\tB\|_F< \|\tA P-P\tB\|_F\};\\
\mathcal{B}_2&:=\{\exists Q\in \mathcal{O}_n^\mathcal{C}\text{ s.t. } \|\tA Q-Q\tB\|_F< \|\tA Q_{\id}-Q_{\id}\tB\|_F\}.
\end{align*}
$\mathcal{B}_1$ is the event that the optimal GMP permutation is not in $\pnc$ and the graphs are not 
core $\delta$-matchable
for 
$\delta(A,PBP^T)=\|A-\e(A)-P(B-\e(B))P^T\|_F$.
As $Q_{\id}\in\pnc$, we have that $\mathcal{B}_1\subset \mathcal{B}_2$.

Suppose that $Q\in\mathcal{O}_n^\mathcal{C}$ (with corresponding permutation $\tau\in S(n)$) permutes $k_c+k_j>0$ core labels, where 
\begin{align}
\label{eq:ccerr}
k_c&=|\{i\in \mathcal{C}: i\neq\tau(i)\in \mathcal{C}\}|,\\
\label{eq:cjerr}
k_j&=|\{i\in \mathcal{C}: i\neq\tau(i)\in \mathcal{J}\}|\\
&=|\{i\in \mathcal{J}: i\neq\tau(i)\in \mathcal{C}\}|.\notag
\end{align}
Applying the results in Appendix \ref{sub:bibern} on the Bivariate Bernoulli distribution, we see that $\|\tA Q-Q\tB\|^2_F- \|\tA Q_{\id}-Q_{\id}\tB\|^2_F$ is a function of $N_P$ independent Bernoulli random variables, where
\begin{align*}
N_P& =3\underbrace{\left(\binom{k_c+k_j}{2}+(n_c-k_c-k_j)(k_c+k_j)\right)}_{\text{core--to--core edges permuted by }\tau\text{ and not }\tau_{\text{id}}} +2\underbrace{\left((k_c+k_j)n_j+(n_c-k_c-k_j)k_j\right)}_{\text{core--to--junk edges permuted differently by }\tau\text{ and }\tau_{\text{id}}}\\
&\hspace{10mm}+2\underbrace{\left(\binom{k_j}{2}+(n_j-k_j)k_j  \right)}_{\text{junk--to--junk edges permuted differently by }\tau\text{ and }\tau_{\text{id}}} \\
&\leq 3(k_c+2k_j)n.
\end{align*}

\tb{As in the proof of Theorem \ref{thm:centeredmatched}, we next apply Proposition \ref{prop:kim} to bound the probability that $Q$ provides a better matching than $Q_{\text{id}}$.
By the assumption that $R(u,v)=0$ if either $u$ or $v$ is a junk vertex, it holds that 
$$\frac{1}{2}\left(\Ex\left(\|\tA Q-Q\tB\|^2_F- \|\tA Q_{\id}-Q_{\id}\tB\|^2_F \right)\right)=\frac{1}{2}\left(\Ex\left(\|\tA Q-Q\tB\|^2_F- \|\tA -\tB\|^2_F \right)\right)=:\mathcal{X}_Q.$$
To ease notation, we define $X_Q:=\frac{1}{2}\left(\|\tA Q-Q\tB\|^2_F- \|\tA Q_{\id}-Q_{\id}\tB\|^2_F\right)$, so that 
\begin{align*}
\Pr\left(X_Q\leq 0\right)&\leq\Pr\left(|X_Q-\mathcal{X}_Q|\geq \mathcal{X}_Q\right)\\
&\leq 2\exp\left\{-\frac{2\mathcal{X}_Q^2}{8N_P}\right\}\\
&\leq 2\exp\left\{ -\frac{2\mathcal{X}_Q^2}{24(k_c+2k_j)n} \right\}.
\end{align*}}

\tb{To use a union bound, note that the number of permutations $Q\in \mathcal{O}_n^\mathcal{C}$ with error counts in Eqs. (\ref{eq:ccerr})--(\ref{eq:cjerr}) given by $k_c$ and $k_j$ is bounded above by $n_c^{2(k_c+k_j)}n_j^{2n_j}$.
Let $k:=k_c+k_j$ be the number of core vertices permuted by $Q$.
Hence, if $n_c\geq n_j$ and
\begin{align*}
\mathcal{X}_Q&=\omega\left( k\sqrt{n_c\log{n_c}}\right),\\
\mathcal{X}_Q&=\omega\left( \sqrt{kn_cn_j\log{n_j}}\right)
\end{align*} then we have that
\begin{align*}
\Pr(\mathcal{B}_1)\leq&\,\,\Pr[ \min_{Q\in \mathcal{O}_n^\mathcal{C}} X_Q \leq 0 ]\\
\leq& \sum_{k_c=1}^{n_c} \sum_{k_j=1}^{\min(n_c-k_c,n_j)} 2 n_c^{2(k_c+k_j)}n_j^{2n_j}\exp\left\{ -\frac{2\mathcal{X}_Q^2}{24(k_c+2k_j)n} \right\}\\
\leq&\sum_{k_c=1}^{n_c} \sum_{k_j=1}^{n_j} 2 \exp \bigg\{-\frac{2\mathcal{X}_Q^2}{48(k_c+k_j)n}+2(k_c+k_j)\log(n_c) +2n_j\log(n_j)\bigg\}\\
\leq& \sum_{k_c=1}^{n_c} \sum_{k_j=1}^{n_j}  2 \exp \left\{-\omega(\log n_c)\right\}
= 2\exp \left\{-\omega(\log n_c)\right\},
\end{align*}
as desired.}


\begin{thebibliography}{10}

\bibitem{Abbe2015-ua}
E.~Abbe and C.~Sandon.
\newblock Community detection in general stochastic block models: Fundamental
  limits and efficient algorithms for recovery.
\newblock In {\em 2015 {IEEE} 56th Annual Symposium on Foundations of Computer
  Science}, pages 670--688, October 2015.

\bibitem{asta2014geometric}
D.~Asta and C.~Shalizi.
\newblock Geometric network comparisons.
\newblock In {\em Proceedings of the Thirty-First Conference on Uncertainty in
  Artificial Intelligence}, UAI'15, pages 102--110, Arlington, Virginia, United
  States, 2015. AUAI Press.

\bibitem{babai2016graph}
L.~Babai.
\newblock Graph isomorphism in quasipolynomial time.
\newblock In {\em Proceedings of the forty-eighth annual ACM symposium on
  Theory of Computing}, pages 684--697. ACM, 2016.

\bibitem{eralg}
B.~Barak, C.~Chou, Z.~Lei, T.~Schramm, and Y.~Sheng.
\newblock (nearly) efficient algorithms for the graph matching problem on
  correlated random graphs.
\newblock {\em arXiv preprint arXiv:1805.02349}, 2018.

\bibitem{netalign}
M.~Bayati, M.~Gerritsen, D.~F. Gleich, A.~Saberi, and Y.~Wang.
\newblock Algorithms for large, sparse network alignment problems.
\newblock In {\em 2009 Ninth IEEE International Conference on Data Mining},
  pages 705--710. IEEE, 2009.

\bibitem{boll}
B.~Bollob{\'a}s.
\newblock {\em Random Graphs}.
\newblock Springer, 1998.

\bibitem{qapref}
S.~Bougleux, L.~Brun, V.~Carletti, P.~Foggia, B.~Ga{\"u}z{\`e}re, and M.~Vento.
\newblock Graph edit distance as a quadratic assignment problem.
\newblock {\em Pattern Recognition Letters}, 87:38--46, 2017.

\bibitem{cela2013quadratic}
E.~Cela.
\newblock {\em The quadratic assignment problem: theory and algorithms},
  volume~1.
\newblock Springer Science \& Business Media, 2013.

\bibitem{chatterjee2015matrix}
S.~Chatterjee.
\newblock Matrix estimation by universal singular value thresholding.
\newblock {\em The Annals of Statistics}, 43(1):177--214, 2015.

\bibitem{chen2016joint}
L.~Chen, J.~T. Vogelstein, V.~Lyzinski, and C.~E. Priebe.
\newblock A joint graph inference case study: the c. elegans chemical and
  electrical connectomes.
\newblock In {\em Worm}, volume~5. Taylor \& Francis, 2016.

\bibitem{jointLi}
L.~Chen, J.T. Vogelstein, V.~Lyzinski, and C.~Priebe.
\newblock A joint graph inference case study: the c.elegans chemical and
  electrical connectomes.
\newblock {\em arXiv preprint, submitted}, 2015.

\bibitem{fcci}
F.~Chung and L.~Lu.
\newblock Concentration inequalities and martingale inequalities: a survey.
\newblock {\em Internet Mathematics}, 3(1):79--127, 2006.

\bibitem{ConteReview}
D.~Conte, P.~Foggia, C.~Sansone, and M.~Vento.
\newblock Thirty years of graph matching in pattern recognition.
\newblock {\em International Journal of Pattern Recognition and Artificial
  Intelligence}, 18(03):265--298, 2004.

\bibitem{het2}
D.~Cullina and N.~Kiyavash.
\newblock Improved achievability and converse bounds for erdos-r{\'e}nyi graph
  matching.
\newblock In {\em ACM SIGMETRICS Performance Evaluation Review}, volume~44,
  pages 63--72. ACM, 2016.

\bibitem{het1}
D.~Cullina and N.~Kiyavash.
\newblock Exact alignment recovery for correlated erd$\backslash$h $\{$o$\}$
  sr$\backslash$'enyi graphs.
\newblock {\em arXiv preprint arXiv:1711.06783}, 2017.

\bibitem{het3}
D.~Cullina, N.~Kiyavash, P.~Mittal, and H.~V. Poor.
\newblock Partial recovery of erd$\backslash$h $\{$o$\}$
  sr$\backslash$'$\{$e$\}$ nyi graph alignment via $ k $-core alignment.
\newblock {\em arXiv preprint arXiv:1809.03553}, 2018.

\bibitem{desikan2006automated}
R.~S. Desikan, F.~S{\'e}gonne, B.~Fischl, B.~T. Quinn, B.~C. Dickerson,
  D.~Blacker, R.~L. Buckner, A.~M. Dale, R.~P. Maguire, and B.~T. Hyman.
\newblock An automated labeling system for subdividing the human cerebral
  cortex on mri scans into gyral based regions of interest.
\newblock {\em Neuroimage}, 31(3):968--980, 2006.

\bibitem{eralg2}
J.~Ding, Z.~Ma, Y.~Wu, and J.~Xu.
\newblock Efficient random graph matching via degree profiles.
\newblock {\em arXiv preprint arXiv:1811.07821}, 2018.

\bibitem{prot}
A.~Elmsallati, C.~Clark, and J.~Kalita.
\newblock Global alignment of protein-protein interaction networks: A survey.
\newblock {\em IEEE/ACM transactions on computational biology and
  bioinformatics}, 13(4):689--705, 2016.

\bibitem{Emmert-Streib2016-st}
F.~Emmert-Streib, M.~Dehmer, and Y.~Shi.
\newblock Fifty years of graph matching, network alignment and network
  comparison.
\newblock {\em Information sciences}, 346--347:180--197, 2016.

\bibitem{escolano}
F.~Escolano, E.~R. Hancock, and M.~Lozano.
\newblock Graph matching through entropic manifold alignment.
\newblock In {\em Computer Vision and Pattern Recognition (CVPR), 2011 IEEE
  Conference on}, pages 2417--2424. IEEE, 2011.

\bibitem{fang2018tractable}
F.~Fang, D.~L. Sussman, and V.~Lyzinski.
\newblock Tractable graph matching via soft seeding.
\newblock {\em arXiv preprint arXiv:1807.09299}, 2018.

\bibitem{Feldman2013-we}
V.~Feldman, E.~Grigorescu, L.~Reyzin, S.~Vempala, and Y.~Xiao.
\newblock Statistical algorithms and a lower bound for detecting planted
  cliques.
\newblock In {\em Proceedings of the Forty-fifth Annual {ACM} Symposium on
  Theory of Computing}, STOC '13, pages 655--664, New York, NY, USA, 2013. ACM.

\bibitem{foggia2014graph}
P.~Foggia, G.~Percannella, and M.~Vento.
\newblock Graph matching and learning in pattern recognition in the last 10
  years.
\newblock {\em International Journal of Pattern Recognition and Artificial
  Intelligence}, 28(01):1450001, 2014.

\bibitem{FW}
M.~Frank and P.~Wolfe.
\newblock An algorithm for quadratic programming.
\newblock {\em Naval Research Logistics Quarterly}, 3(1-2):95--110, 1956.

\bibitem{regal}
M.~Heimann, H.~Shen, T.~Safavi, and D.~Koutra.
\newblock Regal: Representation learning-based graph alignment.
\newblock In {\em Proceedings of the 27th ACM International Conference on
  Information and Knowledge Management}, pages 117--126. ACM, 2018.

\bibitem{hoff2002latent}
P.~D. Hoff, A.~E. Raftery, and M.~S. Handcock.
\newblock Latent space approaches to social network analysis.
\newblock {\em Journal of the american Statistical association},
  97(460):1090--1098, 2002.

\bibitem{sbm}
P.~W. Holland, K.~B. Laskey, and S.~Leinhardt.
\newblock Stochastic blockmodels: First steps.
\newblock {\em Social networks}, 5(2):109--137, 1983.

\bibitem{horn2012matrix}
R.~A. Horn and C.~R. Johnson.
\newblock {\em Matrix analysis}.
\newblock Cambridge university press, 2012.

\bibitem{kazemi2015can}
E.~Kazemi, L.~Yartseva, and M.~Grossglauser.
\newblock When can two unlabeled networks be aligned under partial overlap?
\newblock In {\em Communication, Control, and Computing (Allerton), 2015 53rd
  Annual Allerton Conference on}, pages 33--42. IEEE, 2015.

\bibitem{klau}
G.~W. Klau.
\newblock A new graph-based method for pairwise global network alignment.
\newblock {\em BMC bioinformatics}, 10(1):S59, 2009.

\bibitem{landman2011multi}
B.~A. Landman, A.~J. Huang, A.~Gifford, D.~S. Vikram, I.~A.~L. Lim, J.~A.~D.
  Farrell, J.~A. Bogovic, J.~Hua, M.~Chen, and S.~Jarso.
\newblock Multi-parametric neuroimaging reproducibility: a 3-t resource study.
\newblock {\em Neuroimage}, 54(4):2854--2866, 2011.

\bibitem{le2017concentration}
C.~M. Le, E.~Levina, and R.~Vershynin.
\newblock Concentration and regularization of random graphs.
\newblock {\em Random Structures \& Algorithms}, 51(3):538--561, 2017.

\bibitem{lee2010graph}
J.~Lee, M.~Cho, and K.~M. Lee.
\newblock A graph matching algorithm using data-driven markov chain monte carlo
  sampling.
\newblock In {\em Pattern Recognition (ICPR), 2010 20th International
  Conference on}, pages 2816--2819. IEEE, 2010.

\bibitem{limatching}
L.~Li and W.~M. Campbell.
\newblock Matching community structure across online social networks.
\newblock {\em NIPS Workshop on Networks in the Social and Information
  Sciences}, 2015.

\bibitem{lin2010layered}
L.~Lin, X.~Liu, and S.-C. Zhu.
\newblock Layered graph matching with composite cluster sampling.
\newblock {\em IEEE Transactions on Pattern Analysis and Machine Intelligence},
  32(8):1426--1442, 2010.

\bibitem{loiola2007survey}
E.~M. Loiola, N.~M.~M. de~Abreu, P.~O. Boaventura-Netto, P.~Hahn, and
  T.~Querido.
\newblock A survey for the quadratic assignment problem.
\newblock {\em European journal of operational research}, 176(2):657--690,
  2007.

\bibitem{lyzinski2016information}
V.~Lyzinski.
\newblock Information recovery in shuffled graphs via graph matching.
\newblock {\em IEEE Transactions on Information Theory}, DOI:
  10.1109/TIT.2018.2808999, 2018.

\bibitem{rel}
V.~Lyzinski, D.~E. Fishkind, M.~Fiori, J.~T. Vogelstein, C.~E. Priebe, and
  G.~Sapiro.
\newblock Graph matching: Relax at your own risk.
\newblock {\em IEEE Transactions on Pattern Analysis and Machine Intelligence},
  38(1):60--73, 2016.

\bibitem{JMLR:v15:lyzinski14a}
V.~Lyzinski, D.~E. Fishkind, and C.~E. Priebe.
\newblock Seeded graph matching for correlated {E}rdos-{R}enyi graphs.
\newblock {\em Journal of Machine Learning Research}, 15:3513--3540, 2014.

\bibitem{vncon}
V.~Lyzinski, K.~Levin, and C.~E. Priebe.
\newblock On consistent vertex nomination schemes.
\newblock {\em arXiv preprint arXiv:1711.05610}, 2017.

\bibitem{lyzinski2017consistent}
V.~Lyzinski, K.~Levin, and C.~E. Priebe.
\newblock On consistent vertex nomination schemes.
\newblock {\em Journal of Machine Learning Research}, accepted for publiction,
  2019.

\bibitem{mcdiarmid1989method}
C.~McDiarmid.
\newblock On the method of bounded differences.
\newblock {\em Surveys in combinatorics}, 141(1):148--188, 1989.

\bibitem{Mossel2014-lk}
E.~Mossel, J.~Neeman, and A.~Sly.
\newblock Belief propagation, robust reconstruction and optimal recovery of
  block models.
\newblock In {\em Conference on Learning Theory}, pages 356--370. jmlr.org,
  29~May 2014.

\bibitem{onaran2016optimal}
E.~Onaran, S.~Garg, and E.~Erkip.
\newblock Optimal de-anonymization in random graphs with community structure.
\newblock {\em arXiv preprint arXiv:1602.01409}, 2016.

\bibitem{pedarsani2011privacy}
P.~Pedarsani and M.~Grossglauser.
\newblock On the privacy of anonymized networks.
\newblock In {\em Proceedings of the 17th ACM SIGKDD international conference
  on Knowledge discovery and data mining}, pages 1235--1243. ACM, 2011.

\bibitem{robles}
A.~Robles-Kelly and E.~R. Hancock.
\newblock A riemannian approach to graph embedding.
\newblock {\em Pattern Recognition}, 40(3):1042--1056, 2007.

\bibitem{sang2012robust}
J.~Sang and C.~Xu.
\newblock Robust face-name graph matching for movie character identification.
\newblock {\em IEEE Transactions on Multimedia}, 14(3):586--596, 2012.

\bibitem{sussman2018matched}
D.~L. Sussman, V.~Lyzinski, Y.~Park, and C.~E. Priebe.
\newblock Matched filters for noisy induced subgraph detection.
\newblock {\em arXiv preprint arXiv:1803.02423}, 2018.

\bibitem{MT2}
M.~Tang, A.~Athreya, D.~L. Sussman, V.~Lyzinski, Y.~Park, and C.~E. Priebe.
\newblock A semiparametric two-sample hypothesis testing problem for random dot
  product graphs.
\newblock {\em Journal of Computational and Graphical Statistics},
  26(2):344--354, 2017.

\bibitem{trosset2008semisupervised}
M.~W. Trosset, C.~E. Priebe, Y.~Park, and M.~I. Miller.
\newblock Semisupervised learning from dissimilarity data.
\newblock {\em Computational statistics \& data analysis}, 52(10):4643--4657,
  2008.

\bibitem{udell2017nice}
M.~Udell and A.~Townsend.
\newblock Nice latent variable models have log-rank.
\newblock {\em arXiv preprint arXiv:1705.07474}, 2017.

\bibitem{umeyama1988eigendecomposition}
S.~Umeyama.
\newblock An eigendecomposition approach to weighted graph matching problems.
\newblock {\em IEEE transactions on pattern analysis and machine intelligence},
  10(5):695--703, 1988.

\bibitem{FAQ}
J.~T. {Vogelstein}, J.~M. {Conroy}, V.~{Lyzinski}, L.~J. {Podrazik}, S.~G.
  {Kratzer}, E.~T. {Harley}, D.~E. {Fishkind}, R.~J. {Vogelstein}, and C.~E.
  {Priebe}.
\newblock {Fast Approximate Quadratic Programming for Graph Matching}.
\newblock {\em PLoS ONE}, 10(04), 2014.

\bibitem{xu2017rates}
J.~Xu.
\newblock Rates of convergence of spectral methods for graphon estimation.
\newblock {\em arXiv preprint arXiv:1709.03183}, 2017.

\bibitem{yartseva2013performance}
L.~Yartseva and M.~Grossglauser.
\newblock On the performance of percolation graph matching.
\newblock In {\em Proceedings of the first ACM conference on Online social
  networks}, pages 119--130. ACM, 2013.

\bibitem{young2007random}
S.~Young and E.~Scheinerman.
\newblock Random dot product graph models for social networks.
\newblock In {\em Proceedings of the 5th international conference on algorithms
  and models for the web-graph}, pages 138--149, 2007.

\bibitem{zaslavskiy2009}
M.~Zaslavskiy, F.~Bach, and J.P. Vert.
\newblock A path following algorithm for the graph matching problem.
\newblock {\em Pattern Analysis and Machine Intelligence, IEEE Transactions
  on}, 31(12):2227--2242, 2009.

\bibitem{zhang2016final}
S.~Zhang and H.~Tong.
\newblock Final: Fast attributed network alignment.
\newblock In {\em Proceedings of the 22nd ACM SIGKDD International Conference
  on Knowledge Discovery and Data Mining}, pages 1345--1354. ACM, 2016.

\bibitem{zhang2018consistent}
Y.~Zhang.
\newblock Consistent polynomial-time unseeded graph matching for lipschitz
  graphons.
\newblock {\em arXiv preprint arXiv:1807.11027}, 2018.

\bibitem{zhang2018unseeded}
Y.~Zhang.
\newblock Unseeded low-rank graph matching by transform-based unsupervised
  point registration.
\newblock {\em arXiv preprint arXiv:1807.04680}, 2018.

\bibitem{zhou}
F.~Zhou and F.~De~la Torre.
\newblock Factorized graph matching.
\newblock In {\em Computer Vision and Pattern Recognition (CVPR), 2012 IEEE
  Conference on}, pages 127--134. IEEE, 2012.

\bibitem{zhu2006automatic}
M.~Zhu and A.~Ghodsi.
\newblock Automatic dimensionality selection from the scree plot via the use of
  profile likelihood.
\newblock {\em Computational Statistics \& Data Analysis}, 51(2):918--930,
  2006.

\end{thebibliography}
\end{document}